\documentstyle[11pt]{article}
\input amssym.def
\input amssym
\input epsf

\newtheorem{theorem}{Theorem}[section]
\newtheorem{lemma}[theorem]{Lemma}
\newtheorem{proposition}[theorem]{Proposition}
\newtheorem{corollary}[theorem]{Corollary}
\newtheorem{conjecture}[theorem]{Conjecture}
\newtheorem{definition}[theorem]{Definition}

\newtheorem{problem}[theorem]{Problem}

\begin{document}
\newcommand{\Z}{\Bbb Z}
\newcommand{\R}{\Bbb R} 
\newcommand{\Q}{\Bbb Q}
\newcommand{\C}{\Bbb C}
\newcommand{\lra}{\longrightarrow}
\newcommand{\hra}{\hookrightarrow}
\newcommand{\lms}{\longmapsto}

\begin{titlepage}
\title{Euler complexes   and geometry of modular varieties}
\author{A.B. Goncharov }
\date{{\it To Iosif Bernstein for his 60th birthday}}
\end{titlepage}
\stepcounter{page}
\maketitle

\tableofcontents
\section{Introduction}

\subsection{Summary}
In \cite{G2}, \cite{G3}, \cite{G5}  we described a mysterious connection between the 
depth $m$ 
multiple polylogarithms at $N$-th roots of unity
\begin{equation} \label{Lim}
{\rm Li}_{n_1, ..., n_m}(z_1, ..., z_m):= \sum_{0< k_1 < ... < k_m}
\frac{z_1^{k_1}... z_m^{k_m}}{k^{n_1}_1...k_m^{n_m}}, \quad z_i^N =1,
\end{equation}
and the modular variety 
\begin{equation} \label{Lim1}
Y_1(m,N):= \Gamma_1(m,N)\backslash SL_m(\R)/SO_m
\end{equation}
for the congruence subgroup $\Gamma_1(m;N)$ of $GL_m(\Z)$ 
stabilizing the row $(0, ..., 0,1)$ modulo $N$. 

 The multiple polylogarithms at $N$-th roots of unity provide us all periods of the 
mixed Hodge structure on the pronilpotent completion of the fundamental group 
$\pi_1({\Bbb G}_m - \mu_N, v_0)$, where $\mu_N$ is the group of all $N$-th roots of unity and 
$v_0 = \partial/\partial t$ is the standard tangent vector at zero.

The refined version of the above connection describes the structure of the 
motivic fundamental group 
$\pi_1^{\cal M}({\Bbb G}_m - \mu_N, v_0)$. 
In particular, in the $l$-adic realization it relates the image of the Galois group acting on the pro-$l$ fundamental group 
$\pi_1^{(l)}({\Bbb G}_m - \mu_N, v_0)$  
with the  geometry of the modular varieties (\ref{Lim1}), for all $m$. 

\vskip 3mm

In this paper we give an explanation of the story for $m=2$. 
Recall that the {\it modular complex} for 
$GL_2(\Z)$ is the chain complex of the modular triangulation of the hyperbolic plane 
 ${\cal H}_2$, see Fig. \ref{modd}, placed in degrees $[1,2]$: 
\begin{figure}[ht]
\centerline{\epsfbox{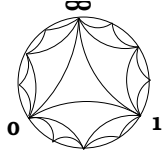}}
\caption{The modular triangulation of the hyperbolic plane.}
\label{modd}
\end{figure}

We assume, for simplicity, that $n_1=n_2 =1$, i.e. work 
with the double logarithm ${\rm Li}_{1,1}(z_1,z_2)$. 
The structure of the (motivic) double 
logarithm at $N$-th roots of unity was described in \cite{G2} 
by a length two complex, called the {\it level $N$ cyclotomic complex}. 
 The above connection
in the  double logarithm case is described by 
a surjective homomorphism 
\begin{equation} \label{al1ui}
\mbox{The modular complex  for $GL_2(\Z)$} \lra 
\mbox{The  level $N$ cyclotomic complex}.
\end{equation}
It is factorized via the coinvariants of the action of the group 
$\Gamma_1(N)$ on the modular complex. 

For every modular curve, we define  
an {\it Euler complex datum}. It is given by 
a map of complexes of vector spaces
\begin{equation} \label{al1}
\mbox{The modular complex for $GL_2(\Z)$} \lra \mbox{The weight two motivic complex on the modular curve.}
\end{equation}
Its image  is called the 
{\it Euler complex} on the modular curve. 
Passing to its second 
 cohomology we recover the 
Beilinson-Kato Euler system in $K_2$ of the tower of modular 
curves (\cite{B}, \cite{Ka}).

We show that the map (\ref{al1ui}) 
can be obtained as the composition of the 
map (\ref{al1}) on the modular curve $Y_1(N)$ 
followed by the specialization at a cusp.

\vskip3mm
At the same time Euler complexes suggest a new twist of the story.
Let us restrict the Euler complex 
to the point corresponding to a CM elliptic curve $E_K$ with complex multiplication by the ring of integers ${\cal O}_K$ in an imaginary quadratic field $K$, and pass to a subcomplex corresponding to the ${\cal N}$-torsion points, where 
${\cal N}$ is an ideal of ${\cal O}_K$. We relate 
the complex obtained to the geometry of the 
three dimensional modular hyperbolic orbifold 
\begin{equation} \label{7.7.05.100}
{\cal Y}_1({\cal N}) := \Gamma_1({\cal N})\backslash SL_2(\C)/SU(2)
\end{equation}
where $\Gamma_1({\cal N}) \subset GL_2({\cal O}_K)$ 
is the subgroup of matrices stabilising the row $(0,1)$ modulo ${\cal N}$. 
The relationship is very precise when $K$ is the field of Gaussian or 
Eisenstein numbers, and the ideal ${\cal N}$ is prime. 
However it is still rather elusive for other imaginary quadratic fields.

In a sequel to this paper we will show that the obtained complex  
is closely related to the $\Q_l(2)$-part of the image of Galois group acting on 
the pro-$l$ completion of the 
 fundamental group of $E_K - \{\mbox{${\cal N}$-torsion points}\}$. 
Thus we get new examples  of the mysterious connection between Galois groups 
 and geometry of modular varieties.
\vskip 3mm
The two examples discussed in the paper can be seen as higher analogs of the 
theory of cyclotomic/elliptic units -- 
we discuss this analogy in the final part of the Introduction.

Finally, there is a similar story in the depth two 
and arbitrary weight situation. Since in this paper we want to present the picture in the simplest form, 
we will elaborate it 
 elsewhere.

\subsection{The double logarithm at roots of unity and modular curves (\cite{G2})} 
The simplest way to state this connection goes via
the dilogarithm story, which we recall now.  
 
Let $F$ be an arbitrary field and $\Z[F^*-\{1\}]$  the free abelian group generated by elements $\{x\}$, where $x\in F^*-\{1\}$. Then there is a version of the Bloch complex 
 \begin{equation} \label{Li2}
\delta_2: {\cal B}_2(F) \lra \Lambda^2 F^*, \qquad 
\delta_2\{x\} = (1-x) \wedge x
 \end{equation}
where $
{\cal B}_2(F)
$ is a certain quotient of  $\Z[F^*-\{1\}]$, see Section 2.1. 
 Recall the dilogarithm function
$$
{\rm Li}_2(z) = -\int_0^z\log(1-t)\frac{dt}{t}.
$$
It has a  single-valued version 
$
{\cal L}_2(z):= {\rm Im} {\rm Li}_2(z) + \arg(1-z)\cdot \log|z|
$, 
providing a homomorphism 
$$
{\cal B}_2(\C) \lra \R, \quad \{z\} \lms {\cal L}_2(z).
$$
We will put the Bloch group in a wider conceptual framework in Section 1.3.

\vskip 3mm
The {\it motivic double logarithms} at $N$-th 
roots of unity were defined (\cite{G5}) as elements
\begin{equation} \label{6.14.05.5}
{\rm Li}^{\cal M}_{1,1}(a,b) 
\in B_2(\Q(\zeta_N))\otimes \Q, \qquad a^N = b^N =1, \quad \zeta_N = 
e^{\frac{2\pi i}{N}}.
\end{equation}
They are motivic avatars of the 
numbers ${\rm Li}_{1,1}(a,b)$. Setting $c:= (ab)^{-1}$, one has 
\begin{equation} \label{6.14.05.5qaz}
\delta_2: {\rm Li}^{\cal M}_{1,1}(a,b) \lms 
(1-a) \wedge (1-b) + (1-b) \wedge (1-c) + (1-c) \wedge (1-a).   
\end{equation}
Therefore 
there is a complex
\begin{equation} \label{Li3}
\delta: C_2(N) \lra \Lambda^2C_1(N).
\end{equation}
where  $C_2(N)$ is the subspace spanned by the elements 
(\ref{6.14.05.5}), and $C_1(N)$ is the 
$\Q$-subspace in 
$\Q(\zeta_N)^*\otimes \Q$ generated by the 
cyclotomic $N$-units $(1-\zeta_N^{\alpha})$. We call it the {\it cyclotomic complex}.
\vskip 3mm
To explain the connection with modular curves, recall 
the modular complex for 
$GL_2(\Z)$:
\begin{equation} \label{4-12.1}
{\rm M}_{(2)}^{\ast}:= {\rm M}_{(2)}^{1} \lra {\rm M}_{(2)}^{2}.
 \end{equation} 
Here ${\rm M}_{(2)}^{1}$ is the group generated by the oriented triangles, 
and ${\rm M}_{(2)}^{2}$ is generated by the oriented geodesics of the modular 
triangulation of the hyperbolic plane shown on Fig \ref{modd}. 
It is a complex of $GL_2(\Z)$-modules. 
Let $\Gamma$ be a subgroup of $GL_2(\Z)$. Projecting the modular 
complex onto the modular curve $Y_{\Gamma}:=  {\cal H}_2/\Gamma$, we get a triangulation of the latter. 
Its chain complex is identified with 
${\rm M}_{(2)}^{\ast}\otimes_\Gamma\Q$. 

\vskip 3mm
Let $\widehat C_1(N) = C_1(N)\oplus \Q$.  Using the embedding $\Lambda^2C_1(N) \hra \Lambda^2\widehat C_1(N)$, we 
 extend the cyclotomic complex (\ref{Li3}) to a complex
\begin{equation} \label{4-12.a1q}
C_2(N) \lra \Lambda^2\widehat C_1(N).
\end{equation} 
We defined in \cite{G2} a canonical 
map from the modular complex to this complex: 
\begin{equation} \label{4-12.a1}
\begin{array}{ccc}
\Bigl({\rm M}_{(2)}^{1}& \lra &{\rm M}_{(2)}^{2}\Bigr)\otimes_{\Gamma_1(N)}\Q\\
\downarrow && \downarrow\\ 
C_2(N) &\lra &\Lambda^2\widehat C_1(N)
\end{array}
\end{equation} 
It enjoys the following properties. The right vertical 
arrow is surjective. It is an isomorphism for a prime $N$.  
The space $C_2(N)$ contains a subspace 
isomorphic to $K_3(\Q(\zeta_N))_\Q$, where $A_\Q:= A\otimes \Q$: 
it is generated by the values 
of the motivic dilogarithm ${\rm Li}^{\cal M}_2$ 
at $N$-th roots of unity. It is the kernel 
of the bottom arrow. 
Modulo the subspace $K_3(\Q(\zeta_N))_\Q$, the left vertical map 
is surjective for all $N$, and is an isomorphism
 for a prime $N$. Therefore for a prime level $N=p$,  
the  map of complexes (\ref{4-12.a1}) gives rise 
to an isomorphism of the top complex with
 the quotient of the bottom one by the subspace $K_3(\Q(\zeta_N))_\Q$

\vskip 3mm

\subsection{Mixed Tate motives and the Bloch group}
Let $S$ be a set of primes in a number field $F$, and 
${\cal O}_{F,S}$ the ring of $S$-integers in $F$. Then there is an abelian $\Q$-category 
${\cal M}_T({\rm Spec}{\cal O}_{F,S})$ of mixed Tate motives over ${\rm Spec}{\cal O}_{F,S}$, defined in \cite{DG}. 
It is equipped with a canonical fiber functor 
$$
\omega: {\cal M}_T({\rm Spec}{\cal O}_{F,S}) \lra \mbox{$\Q$-vector spaces}.
$$ 
Let ${\rm L}_{\bullet}({\rm Spec}{\cal O}_{F,S})$ be the Lie algebra  of its derivations, 
called the fundamental Lie algebra of ${\cal O}_{F,S}$\footnote{We apologize for the abuse of notation:  ${\rm L}_{2}(z)$ denotes the 
Bloch-Wigner, and both notation are rather standard.} 
 It is  graded by negative integers. The functor $\omega$ provides an equivalence between the category 
of mixed Tate motives over ${\rm Spec}{\cal O}_{F,S}$ and the category 
of graded finite dimensional modules over the fundamental Lie algebra. 
Let us denote by ${\cal L}_{\bullet}({\rm Spec}{\cal O}_{F,S})$ 
the graded dual of the fundamental Lie algebra. 
It is a positively graded Lie coalgebra. We will usually omit ${\rm Spec}$ from notation. 

The fundamental Lie algebra ${\rm L}_{\bullet}({\cal O}_{F,S})$ is a free graded Lie algebra 
generated by $K_{2n-1}({\cal O}_{F,S})_\Q$ in degree $n$, where $n =1, 2, ...$ ({\it loc. cit.}). 
So if the set $S$ is finite all graded components are finite dimensional. 
In particular, setting $n=1$, we arrive at  a  canonical isomorphism
\begin{equation} \label{10.15.05.4}
{\cal O}_{F,S}^*\otimes \Q \stackrel{\sim}{\lra}{\cal L}_{1}({\cal O}_{F,S}). 
\end{equation}
Furthermore, there is a canonical isomorphism (Lemma \ref{10.15.08.10q})
\begin{equation} \label{10.15.05.8}
{\cal B}_2(F)_\Q \stackrel{\sim}{\lra} {\cal L}_{2}(F).
\end{equation}
Let ${\cal B}_2({\cal O}_{F,S})_\Q $ be the biggest subspace in ${\cal B}_2(F)_\Q$ which is mapped by 
the differential $\delta_2$ to $\Lambda^2{\cal O}_{F,S}^*\otimes \Q$. Then 
the isomorphism (\ref{10.15.05.8}) restricts (Lemma \ref{10.15.08.10}) to an isomorphism
\begin{equation} \label{10.15.05.9}
{\cal B}_2({\cal O}_{F,S})_\Q 
\stackrel{\sim}{\lra} {\cal L}_{2}({\cal O}_{F,S}).
\end{equation}
\vskip 3mm
Let us introduce the two cyclotomic schemes
$$
S_N:= {\rm Spec} \Z[\zeta_N][\frac{1}{N}], \qquad S^{\rm un}_N:= {\rm Spec} \Z[\zeta_N].
$$
The cyclotomic complex (\ref{Li3}) is a subcomplex in the weight two part 
$$
{\cal L}_{2}(S_N) \lra \Lambda^2 {\cal L}_{1}(S_N)
$$
of the standard 
cochain complex of the fundamental Lie coalgebra ${\cal L}_{2}(S_N)$. Moreover, if $N=p$ is a prime, 
we can replace  here $S_p$ by $S_p^{\rm un}$. 
Analysing the diagram 
(\ref{4-12.a1}) we arrive at the following 
\vskip 3mm
{\bf Conclusion 1 (\cite{G2}).} {\it For a prime $p$, the subspace of 
 ${\cal L}_{2}(S^{\rm un}_p)$ generated by the motivic double logarithms at $p$-th roots of unity 
is smaller then ${\cal L}_{2}(S^{\rm un}_p)$;  the quotient is isomorphic to 
 $H_{\rm cusp}^1(\Gamma_1(p), \Q)$}.
\vskip 3mm 
 Below we develop an analog of this picture related to a CM elliptic curve, see Section 1.6.

\subsection{Relation with the motivic fundamental group of ${\Bbb G}_m - \mu_N$} 
The unipotent 
motivic fundamental group  $\pi_1^{\cal M}({\Bbb G}_m - \mu_N, v_0)$ 
is not a group but rather a pro-unipotent 
Lie algebra object in the abelian category  
of mixed Tate motives over the scheme 
$S_N$ (\cite{DG}). Thus applying the fiber functor $\omega$ 
to the (pro-) mixed Tate motive $\pi_1^{\cal M}({\Bbb G}_m - \mu_N, v_0)$ 
we get a graded 
Lie algebra $\omega(\pi_1^{\cal M}({\Bbb G}_m - \mu_N, v_0))$ 
over $\Q$. 
The fundamental Lie algebra ${\rm L}_{\bullet}(S_N)$ 
acts  by its derivations, 
providing a canonical Lie algebra homomorphism
$$
{\rm L}_{\bullet}(S_N) \lra {\rm Der}\Bigl(\omega(\pi_1^{\cal M}({\Bbb G}_m - \mu_N, v_0))\Bigr).
$$
Let us denote by 
$C_{\bullet}(\mu_N)$ its image. We called it {\it the cyclotomic Lie algebra}. 
Let $C_{\bullet}(\mu_N)^{\vee}$ be the dual Lie coalgebra. One has 
$$
C_{-1}(\mu_N)^{\vee} = 
(\mbox{the group of cyclotomic units in $S_N$})\otimes \Q = {\cal O}_{S_N}^*\otimes \Q.
$$
We proved in \cite{G5} that 
$$
C_{-2}(\mu_N)^{\vee} = C_2(N)
$$
and the dual to the commutator map $[*,*]: \Lambda^2C_{-1}(\mu_N) \lra C_{-2}(\mu_N)$ 
is identified with the cyclotomic complex (\ref{Li3}). Thus
 the cyclotomic complex describes the weight $2$ part 
of the image of the motivic Galois group acting on $\pi_1^{\cal M}({\Bbb G}_m - \mu_N, v_0)$. 
In general, 
the cyclotomic Lie algebra $C_{\bullet}(\mu_N)$ is   
described by 
the motivic multiple polylogarithms at $N$-th roots of unity.

\subsection{Euler complexes on modular curves} 
\begin{figure}[ht]
\centerline{\epsfbox{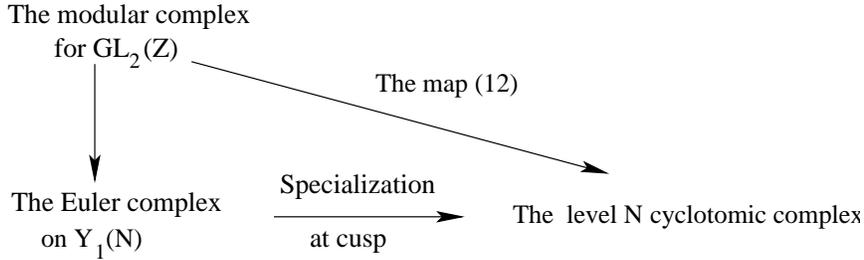}}
\caption{Relating the modular, Euler and cyclotomic complexes.}
\label{ec1}
\end{figure}

In this paper we develop this story as follows. Using 
a construction from Section 6 of \cite{G1}, 
we define 
a modular deformation of the complex (\ref{Li3}), 
called the Euler complex. 
It is a subcomplex of the Bloch complex (\ref{Li2}) for the field of functions 
on the modular curve $Y(N)$. (In fact it is a subcomplex of the 
Bloch complex of the modular curve from Definition 2.1).
We show in Proposition  \ref{6.16.05.1} that there is a canonical map 
\begin{equation} \label{7.7.05.10}
\mbox{the modular complex} \lra \mbox{the Euler complex on $Y(N)$}. 
\end{equation}
One gets a similar map for other modular curves, e.g. $Y_1(N)$. Taking its  
specialization at a cusp $\infty$ 
on $Y_1(N)$ obtained by projection of the point 
$\infty$ from the hyperbolic plane, we recover 
the cyclotomic complex (\ref{Li3}): 
the specialization  provides 
 a surjective morphism of complexes:
\begin{equation} \label{7.7.05.11}
\mbox{The Euler complex on $Y_1(N)$} \stackrel{}{\lra} 
\mbox{The level $N$ cyclotomic complex (\ref{Li3})}.
\end{equation}
This map intertwines the maps from the modular complex to the Euler and cyclotomic complexes. 
So we arrive at  
 the commutative diagram on Fig. \ref{ec1}. 
This  
explains why the modular curve $Y_1(N)$ appears in the study of 
the motivic double logarithm  at roots of unity. 

\vskip 3mm
{\it Generalizations}. Recall the standard cochain complex
\begin{equation} \label{7.7.06.1}
C_{\bullet}(\mu_N)^{\vee} \lra \Lambda^2C_{\bullet}(\mu_N)^{\vee} \lra 
\Lambda^3C_{\bullet}(\mu_N)^{\vee} \lra \ldots 
\end{equation}
of the Lie algebra $C_{\bullet}(\mu_N)$. The first map is dual to the commutator map, and the others are defined using the Leibniz rule. 
The grading of the Lie algebra $C_{\bullet}(\mu_N)$ provides a weight 
grading of the complex. 
The cyclotomic complex (\ref{Li3}) 
is isomorphic to the weight two part of the standard cochain complex 
(\ref{7.7.06.1}) 
of the cyclotomic Lie algebra $C_{\bullet}(\mu_N)$.  

The Lie algebra $C_{\bullet}(\mu_N)$ has a depth filtration.
 So the complex (\ref{7.7.06.1}) inherits the depth filtration. 
The depth $m$ part of $C_{\bullet}(\mu_N)$ is described by 
motivic multiple polylogarithms of the depth $\leq m$. 

In \cite{G3}, see also \cite{G4}, \cite{G5}, we generalized 
 the 
diagonal arrow in Fig. \ref{ec1} to $GL_m(\Z)$, for any positive integer $m$. 
Namely, we defined the rank $m$ modular complex and constructed 
a map from this complex, tensored by $S^{w-m}V_m$, where $V_m$ 
is the standard representation of $GL_m$, to the depth $m$, weight $w$ part of 
the standard cochain complex of the cyclotomic Lie algebra $C_{\bullet}(\mu_N)$. 

What might 
play the role of the Euler complex for $m>2$? Observe that 
the modular varieties for $m>2$ are not algebraic varieties.

\subsection{Euler complexes, CM points, and 
modular hyperbolic $3$-folds}
Let ${\cal O}_K$ be the ring of integers in the imaginary quadratic field 
$\Q(\sqrt{-d})$, and ${\cal N}$ an ideal in ${\cal O}_K$. 
The group $GL_2({\cal O}_K)$ acts 
discretely on the 
hyperbolic $3$-space ${\cal H}_3$. Recall the subgroup 
$\Gamma_1({\cal N})$  of $GL_2({\cal O}_K)$ 
stabilizing the row vector $(0,1)$ modulo ${\cal N}$, and  the corresponding 
modular orbifold ${\cal Y}_1({\cal N}):= 
\Gamma_1({\cal N})\backslash {\cal H}_3$.  
Let $E_K$ be an elliptic curve with the endomorphism ring ${\cal O}_K$. 

\vskip 3mm
Restricting the Euler complex on 
$Y(N)$, where $N$ is the norm of ${\cal N}$,  to a point corresponding to the curve $E_K$, 
and taking the subcomplex corresponding to the ${\cal N}$-torsion points, 
we get a complex which is 

\begin{itemize} 
\item  Related to the geometry of the three-dimensional modular orbifold 
${\cal Y}_1({\cal N})$. The relationship is 
 most precise 
when the ideal ${\cal N}$ is prime, and $K$ is 
the field of Gaussian or Eisenstein numbers. 

\item Related to the simplest (depth two, weight two) 
non-abelian quotient of the image of the Galois group acting on the 
pro-$l$ fundamental group of $E_K- E_K[{\cal N}]$. 
\end{itemize}

Below we elaborate  the first connection.  The second  will be discussed 
in a sequel to this paper.

\paragraph{The modular and elliptic units.} 
Recall that the classical modular unit    is an invertible 
function $\theta_q(z)$ on the modular curve $Y_1(N)$ constructed as follows.  
Let 
$$
q=e^{2 \pi i \tau}, {\rm Im} (\tau) >0, \qquad z=e^{2 \pi i \xi}, \qquad \xi = \alpha_1 \tau 
+ \alpha_2,\quad \alpha_1,\alpha_2 \in \frac{1}{N}\Z^2 - \Z^2
$$
Set
\begin{equation} \label{eq1}
\theta_q(z) = -q^{\frac{1}{2}B_2(\alpha_1)}\cdot e^{2\pi i \alpha_2(\alpha_1-1)/2} (1-z) 
\prod_{n=1}^{\infty}(1-q^nz)(1-q^nz^{-1})
\end{equation}
where $B_2(x):= x^2- x +\frac{1}{6}$ is the second Bernoulli polynomial. 
Changing $(\alpha_1,\alpha_2)$ by an element of $\Z^2$, we alter $\theta_q(z)$ 
by multiplication by an $N$-th root of unity. So $\theta_q(z)^N$ is a well defined  invertible function on $Y_1(N)$. 
  
Let $E$ be an elliptic curve over an arbitrary field $k$. 
Then, given an $N$-torsion point $z$ on $E$ defined over a field $k_z$, 
one can define an element $\theta_E(z)$ 
 such that  $\theta_E(z)^N \in k_z^*$ (see \cite{GL},
 and Section 2.3 below). More generally, if $E$ is an elliptic 
curve over a base $S$, and $z$ is an $N$-torsion section, 
there is an element $\theta_E(z)$ such that 
$\theta_E(z)^N \in {\cal O}_S^*$. 
There are the following interesting special cases: 

\begin{itemize}

\item
If $E$ is the universal elliptic curve ${\cal E}_N$ over  $Y_1(N)$ 
we get the described above modular unit. 

\item If $E$ is the  CM curve $E_K$, and ${\cal N}$ is an ideal of ${\cal O}_K$ with the norm $N$, the  elements   $\theta_E(z)$ for all  ${\cal N}$-torsion points $z$ 
of $E_K$ generate an abelian extension $K_{{\cal N}}$ of $K$. 
Moreover,   $\theta_E(z)$ is an $N$-unit  in   $K_{{\cal N}}$, 
called an elliptic unit.
\end{itemize}

\paragraph{The elements $\theta_E(a,b,c)$.}
For any triple  of 
torsion points $a,b,c$ on an elliptic curve $E$ over an arbitrary field 
$k$ with $a+b+c=0$ we construct  an element  
$$
\theta_E(a,b,c)
 \in {\cal B}_2(k')\otimes \Q
$$
 where $k'$ is the field generated  over $k$ by the 
coordinates of the points $a,b$. An 
explicit construction of this element is  
given by 
a reciprocity law from Section 6 of  \cite{G1}, which strengthens Suslin's  
reciprocity law for the Milnor $K_3$-group  of the function field of 
$E$. 
One has  the following key formula 
\begin{equation} \label{6.13.05.1}
\delta: \theta_E(a,b,c) \lms  \theta_E(a) \wedge \theta_E(b) + \theta_E(b) \wedge \theta_E(c) +  \theta_E(c) \wedge \theta_E(a).
\end{equation}

\paragraph{Double modular units and Euler complexes.} 
In particular, let ${\cal E}_N$ be the universal elliptic curve over the 
modular curve $Y(N)$. Then a pair of torsion sections $a,b$ of 
${\cal E}_N$ provides an element 
$\theta_{{\cal E}_N}(a,b,c)$. It can be thought of as a double modular unit: 
its coproduct is a sum of the wedge products of the classical modular units  
$\theta_{{\cal E}_N}(a)$. 
We define the Euler complex as the subcomplex of the 
complex (\ref{Li2}) on $Y(N)$ spanned by 
$\theta_{{\cal E}_N}(a,b,c)$ in degree $1$, 
and  $\theta_{{\cal E}_N}(a) \wedge \theta_{{\cal E}_N}(b)$ in degree $2$, 
for all torsion sections $a,b$ as above. 

\vskip 3mm
\paragraph{An analytic version of the double modular units.} 
Let $X$ be a complex algebraic variety 
with the function field $\C(X)$. 
An element $\beta$ of the Bloch group ${\cal B}_2(\C(X))$ gives rise to a multivalued analytic function ${\rm Li}_2(\beta)$ at the generic point of $X(\C)$. 
Namely, if $\beta = \sum_i n_i \{f_i(z)\}_2$, we set 
  ${\rm Li}_2(\beta) = \sum_i n_i {\rm Li}_2(f_i(z))$. 
So the element $\theta_{{\cal E}_N}(a,b,c)$ gives rise to a 
multivalued analytic function $\widetilde \theta_{{\cal E}_N}(a,b,c)$ 
at the generic point of the 
modular curve $Y(N)\otimes_{\Q(\mu_N)}\C$.  It follows from 
(\ref{6.13.05.1}) that its differential is  
$$
d\widetilde \theta_{{\cal E}_N}(a,b,c) = {\rm Cycle}_{a,b,c}\Bigl( \log 
\theta_{{\cal E}_N}(a)
d\log \theta_{{\cal E}_N}(b) - \log \theta_{{\cal E}_N}(b)
d\log \theta_{{\cal E}_N}(a) \Bigr).
$$
So 
 elements $\theta_{{\cal E}_N}(a,b,c)$ allow to 
express the integral of the $1$-form on the 
right via the dilogarithm. 

 \paragraph{Double elliptic units and modular $3$-folds.} 
Let us 
specialize the Euler complex on $Y(N)$  to the point 
corresponding to a CM curve $E_K$ with complex 
multiplication by ${\cal O}_K$, and consider only 
${\cal N}$-torsion points, where ${\cal N}$ is an ideal 
of ${\cal O}_K$. Let us spell the definition of the 
obtained complex. 

\begin{definition} The $\Q$-vector space 
${\cal C}_2({\cal N})$ is the subspace of  
${\cal B}_2(K_{{\cal N}})_\Q$ generated by the elements $\theta_E(a,b,c)$ when $a,b,c$ run through all 
${\cal N}$-torsion points of $E$ such that $a+b+c=0$. 
\end{definition}  

Let us denote by ${\cal C}_1({\cal N})$ the group of ${\cal N}$-units in  
the field $K_{\cal N}$ tensored by $\Q$.  It is known to be generated 
by the elliptic units $\theta_E(a)$, when  $a$ runs through the ${\cal N}$-torsion 
points of the curve $E_K$.

It follows from (\ref{6.13.05.1}) that we get a complex 
\begin{equation} \label{mc11}
{\cal C}_2({\cal N}) \stackrel{\delta}{\lra} 
\Lambda^2 {\cal C}_1({\cal N}).
\end{equation} 

\begin{figure}[ht]
\centerline{\epsfbox{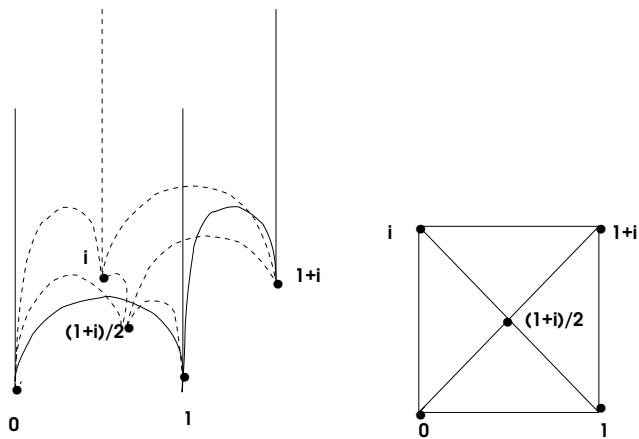}}
\caption{An octahedron of the Bianchi tessellation; 
on the right shown the view from the infinity.}
\label{kj6.00}
\end{figure}

\begin{figure}[ht]
\centerline{\epsfbox{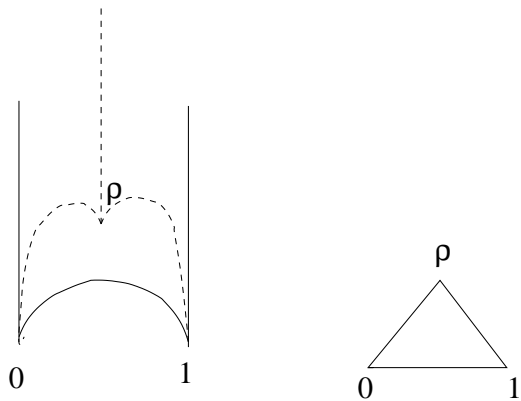}}
\caption{A tetrahedron of the Bianchi tessellation 
for the Eisenstein integers; 
the right picture shows the base viewed from infinity.}
\label{kj7}
\end{figure}

Let ${\cal H}^*_3:= 
{\cal H}_3\cup K \cup \{\infty\}$. We relate this complex with the geometry of the 
modular hyperbolic  orbifold ${\Gamma}_1({\cal N})\backslash {\cal H}^*_3$, 
as follows. 
The hyperbolic space ${\cal H}_3$ has a classical 
Bianchi tessellation on geodesic polyhedrons invariant under 
the action of $GL_2({\cal O}_K)$. For example, for the ring of Gaussian integers   it is a tessellation on octahedron's and for the ring of Eisenstein integers  it is tessellation on tetrahedrons, see Fig. \ref{kj6.00} 
and Fig. \ref{kj7}.

Projecting the Bianchi tessellation onto the modular 
orbifold $\Gamma_1({\cal N})\backslash {\cal H}^*_3$ we 
get a cell decomposition of the latter. 
Let $V_{\bullet}({\cal N})$ be its integral chain complex,  
called the  {\it Bianchi complex} of $\Gamma_1({\cal N})\backslash {\cal H}^*_3$. 
In Section 5 we  relate the complex (\ref{mc11}) 
with the complex $V_{\bullet}({\cal N})$. 
Our results are complete when ${\cal N} = {\cal P} $ 
is a prime ideal,  and $d=-1, -3$. In these cases we 
define a canonical morphism of complexes
\begin{equation} \label{10.16.05.1}
\begin{array}{ccccccc}
V_{3}({\cal P})&\lra &V_{2}({\cal P})&\lra &V_{1}({\cal P})&\lra &V_{0}({\cal P})\\
\downarrow&&\downarrow&&\downarrow&&\downarrow\\
K_3(K_{\cal P})_\Q &\lra &{\cal C}_2({\cal P})&\stackrel{\delta}{\lra} &\Lambda^2 {\cal C}_1({\cal P})
&\stackrel{v}{\lra}&{\cal C}_1({\cal P}) 
\end{array}
\end{equation}
(where $v$ is the residue map related to the valuation defined by ${\cal P}$) and prove that it is almost a  quasiisomorphism -- see 
a precise statement in Theorem \ref{mtheor2da}. 
It is an 
analog of homomorphism (\ref{4-12.a1}). Here is a  corollary. 
Set ${\cal C}_1^{\rm un}({\cal P}):= {\cal O}^*_{K_{\cal P}}\otimes \Q$. Clearly ${\rm Im}\delta \subset \Lambda^2{\cal C}_1^{\rm un}({\cal P})$ in the bottom 
line of (\ref{10.16.05.1}). 

\begin{theorem} \label{10.15.05.2} Let $d=-1$ or $d=-3$, and ${\cal P}$ be a prime ideal in ${\cal O}_K$. Then  
there is an isomorphism
$$
H_{\rm cusp}^2(\Gamma_1({\cal P}), \Q) = {\rm Coker} \Bigl({\cal C}_2
({\cal P})
\lra \Lambda^2 {\cal C}^{\rm un}_1({\cal P}) \Bigr).
$$
\end{theorem}

\paragraph{Double elliptic units and mixed Tate motives over ${\cal O}_{K_{{\cal N}}}$.} 
Set 
\begin{equation} \label{10.15.05.1}
S_{\cal N}:= {\rm 
Spec}{\cal O}_{K_{{\cal N}}}[\frac{1}{N}], \qquad S^{\rm un}_{\cal N}:= {\rm 
Spec}{\cal O}_{K_{{\cal N}}}.
\end{equation} 
Recall the fundamental Lie coalgebra 
${\cal L}_{\bullet}(S_{\cal N})$.  
It follows from (\ref{10.15.05.8}), the very definition of ${\cal C}_2({\cal N})$ as a subspace of 
${\cal B}_2(K_{\cal N})$, 
the key formula (\ref{6.13.05.1}), and Lemma 2.2 below 
that there is an inclusion
$$
i: {\cal C}_2({\cal N}) \hra {\cal L}_{2}(S_{\cal N}).
$$
It gives rise to a homomorphism of complexes
$$
\begin{array}{ccc}
{\cal C}_2({\cal N})&\lra & \Lambda^2{\cal C}_1({\cal N})\\
i\downarrow && \downarrow =\\
{\cal L}_{2}(S_{\cal N}) &\stackrel{}{\lra} & \Lambda^2{\cal L}_{1}(S_{\cal N})
\end{array}
$$
where the right arrow is provided by the isomorphism (\ref{10.15.05.4}) 
for the scheme $S_{\cal N}$, and the bottom arrow dualises the commutator map. 
The kernel of the bottom arrow is identified with $K_3(K_{\cal N})_\Q$. 
Moreover, if ${\cal P}$ is a prime ideal, then we have $i: 
{\cal C}_2({\cal P}) \hra {\cal L}_{2}(S^{\rm un}_{\cal P})$, and  Theorem \ref{10.15.05.2} implies 
\begin{corollary} \label{10.15.05.3} Let $d=-1$ or $d=-3$, and ${\cal P}$ be a prime ideal in ${\cal O}_K$. Then  
there is an isomorphism
\begin{equation}
H_{\rm cusp}^2(\Gamma_1({\cal P}), \Q)\stackrel{\sim}{\lra} \frac{{\cal L}_{2}(S^{\rm un}_{\cal P})}{i({\cal C}_2({\cal P})) 
+ K_3(K_{\cal P})_\Q}. 
\end{equation}
\end{corollary}
Therefore we arrive at the CM analog of the Conclusion 1 in Section 1.3.
\vskip 3mm
{\bf Conclusion 2.} {\it Let $d=-1$ or $d=-3$ and a prime ideal ${\cal P}$ in ${\cal O}_K$. 
Then the subspace of ${\cal L}_2(S^{\rm un}_{\cal P})$ 
generated by the double elliptic units 
at ${\cal P}$-torsion points of $E_K$ 
and $K_3(K_{\cal P})_\Q$ 
is smaller then 
${\cal L}_2(S^{\rm un}_{\cal P})$; the quotient is isomorphic to $H_{\rm cusp}^2(\Gamma_1({\cal P}), \Q)$}. 

\vskip 3mm
{\bf Coda.} The cyclotomic units provide a finite index subgroup in the unit group of a cyclotomic 
field. 
There are two generalizations, highlighted  in Conclusions 1 and 2, where the unit group  is replaced by 
the $\Q(2)$-part of a fundamental Lie coalgebra, 
the role of the cyclotomic unit subgroup is played by an explicitly defined  ``double units'' subspace 
of the fundamental Lie coalgebra, and 
the gap between them is isomorphic to the cuspidal cohomology 
of certain modular varieties. 

In both cases the ``double units'' subspaces can be obtained by specializations of the {\it first} 
groups of the Euler complexes. On the other hand  the 
Beilinson-Kato Euler system is delivered by the {\it second} cohomology groups of the Euler complexes. 
Thus both parts of the Euler complex are important. 

\vskip 3mm 
 {\bf Acknowledgments}. This work was written during my stay at 
IHES in June 2005. 
The key constructions were worked at the MPI(Bonn).
I am grateful to both institutions for hospitality and support. 
I was supported by the NSF grant DMS-0400449.

I am very grateful to Andrey Levin and an anonimous referee, 
who read the 
paper very carefully, made a lot of useful 
remarks, 
 and corrected many misprints/errors. 

\section{Euler complexes on modular curves}

\subsection{The weight two motivic complex and the dilogarithm} 
Let $F$ be an arbitrary field and $\Z[{\Bbb P}^1(F)]$  the free abelian group generated by the elements $\{x\}$ where $x\in {\Bbb P}^1(F)$. Consider a homomorphism $$
\delta_2: \Z[{\Bbb P}^1(F)] \lra \Lambda^2 F^*, \qquad \{x\} \lms (1-x) \wedge x, \quad \{0\}, \{1\}, \{\infty\} \lms 0.
$$ 
  Let ${\cal R}_2(F)$ be the subgroup of $\Z[{\Bbb P}^1(F)]$ 
generated by the following elements. Let $X$ be a  connected curve over $F$ and $f_i$  rational functions on $X$ such that $\sum_i (1-f_i) \wedge  f_i  =0 $ in $\Lambda^2 F(X)^*$. Then ${\cal R}_2(F)$ is generated by  the elements $\sum_i\{f_i(0)\} - \{f_i(1)\}$ for all possible curves  and functions $f_i$ as above. One can show that $\delta_2({\cal R}_2(F))=0$, so setting
$$
{\cal B}_2(F):= \frac{\Z[{\Bbb P}^1(F)]}{{\cal R}_2(F)}
$$
we get a complex (a version of the Bloch complex, see \cite{G6})
$$
\delta_2: {\cal B}_2(F) \lra \Lambda^2 F^*.
 $$

There is a more explicit version $B_2(F)$ 
of the group ${\cal B}_2(F)$. 
Denote by $R_2(F)$  the subgroup of $\Z[{\Bbb P}^1(F)]$ 
generated by $$
\{0\}, \{\infty\}\quad \mbox{and} \quad \sum_{i=1}^{5}(-1)^i
\{r(x_1,...,\widehat   x_i,...,x_5)\},
$$ 
where $(x_1,...,x_5)$ runs through all $5$-tuples of 
distinct points in ${\Bbb P}^1(F)$, where $r(...)$ is the cross-ratio.  The Bloch group $B_2(F)$ is 
the quotient of $\Z[{\Bbb P}^1(F)]$ by the subgroup $R_2(F)$. We get 
the Bloch complex
$$
\delta_2: B_2(F) \lra \Lambda^2F^*.
$$
According to the theorems of Matsumoto and Suslin, one has for this complex:
\begin{equation} \label{10.15.05.9f}
{\rm Coker}\delta_2 = K_2(F), \quad {\rm Ker}\delta_2\otimes \Q = 
K^{\rm ind}_3(F)_\Q .
\end{equation} 
One can show that $R_2(F)\subset {\cal R}_2(F)$. 
Thus there is a map 
$i: B_2(k)\lra {\cal B}_2(k)$ 
induced by the identity map on the generators. 
According Proposition 6.1 of \cite{G1}, this map   is an isomorphism modulo torsion for a 
number field $F$. 
 The function ${\cal L}_2$ provides a homomorphism 
${\cal B}_2(\C) \lra \R$, $\{z\} \lms {\cal L}_2(z)$.

\begin{definition} 
The Bloch group ${\cal B}_2(X)$ of an irreducible  scheme $X$ 
with the field of functions $F_X$ is the largest 
subgroup of ${\cal B}_2(F_X)$ which has the property 
that $\delta_2({\cal B}_2(X)) \subset \Lambda^2{\cal O}_X^*$.
\end{definition}
By the very definition, there is a commutative diagram, where the vertical arrows are embeddings:
$$
\begin{array}{ccc}
{\cal B}_2(X)& \stackrel{\delta_2}{\lra} &\Lambda^2{\cal O}_X^*\\
\downarrow && \downarrow \\
{\cal B}_2(F_X)&\stackrel{\delta_2}{\lra}&\Lambda^2F_X^*
\end{array}
$$

\begin{lemma} \label{10.15.08.10q} Let $F$ be a number field. Then there is a canonical isomorphism
\begin{equation} \label{10.15.05.9}
{\cal B}_2(F)_\Q \stackrel{\sim}{\lra} {\cal L}_{2}(F).
\end{equation}
\end{lemma}

{\bf Proof}. 
Recall the canonical isomorphism ${\cal L}_1(F) = F^*_\Q$. 
Let us define a canonical map $\Z[F] \lra {\cal L}_2(F)$. We 
 assign to a generator $\{x\}$ 
a framed mixed Tate motive ${\rm Li}_2^{\cal M}(x)$ defined as follows. 
Consider ${\Bbb P}^2$ with a projective coordinate system $(z_0, z_1, z_2)$. 
Let  $L_i$ be the line defined by the equation $z_i=0$. 
We get a coordinate triangle $(L_0, L_1, L_2)$. Consider another triple of 
lines $(M_0, M_1, M_2)$, where $M_0:=\{z_1+z_2=1\}$, 
$M_1:=\{z_1=1\}$, $M_2:=\{z_2=x\}$. Let us assume $x \not =0$. 
Let us blow up 
all points where three of the $L$ and $M$-lines intersect: 
the vertices of the $L$-traingle at infinity, the vertex of the 
$M$-triangle on the line $L_2$, and, if there is a vertex of the 
$L$-triangle on the line $L_1$, i.e. $x=1$, this vertex as well.  
Denote by $\widehat {\Bbb P}^2$ the obtained surface. 
Then we get an $M$-pentagon ${\Bbb M}$ 
there, formed by the 
strict preimage of the $M$-triangle plus exceptional lines over the $M$-vertices at infinity. Further, we get $M$-polygon ${\Bbb L}$ (it is a 
pentagon if $x=1$ and $4$-gon otherwise). We define now a 
mixed Tate motive by using the formula 
$$
{\rm Li}_2^{\cal M}(x):= H^2(\widehat {\Bbb P}^2-{\Bbb L}, {\Bbb M} - 
({\Bbb L} \cap {\Bbb M} )).
$$
Its framing, and 
 interpretation of the obtained object 
as a mixed Tate motive is standard, 
see for example \cite{G7}. 
Finally, according to the standard formalism, 
the framed mixed Tate motive ${\rm Li}_2^{\cal M}(x)$ gives rise 
to an element of ${\cal L}_2(F)$.

It is well known, and easy to prove, 
that the coproduct $\delta {\rm Li}_2^{\cal M}(x)$ 
equals $(1-x)\wedge x$. 
So we get a commutative diagram, where the right arrow is an isomorphism 
after $\otimes \Q$:
$$
\begin{array}{ccc}
\Z[F] & \stackrel{\delta_2}{\lra} & \Lambda^2 F^*\\
\downarrow &&\downarrow \\
{\cal L}_2(F)&\stackrel{\delta}{\lra}& \Lambda^2 {\cal L}_1(F)
\end{array}
$$
 Further,  for a complex embedding 
$\sigma: F \hra \C$ the real period of ${\rm Li}_2^{\cal M}(\sigma(x))$ is given by ${\cal L}_2(\sigma(x))$. 
This, the injectivity of the regulator map on $K_3(F)_\Q$, and Suslin's theorem (\ref{10.15.05.9f}) imply
 that the subspace  ${\cal R}_2(F)_\Q \subset \Q[F] $ is killed by the left arrow. 
Thus we get a commutative diagram
$$
\begin{array}{ccc}
{\cal B}_2(F)_\Q & \stackrel{\delta_2}{\lra} & \Lambda^2 F^*\\
\downarrow &&\downarrow \\
{\cal L}_2(F)&\stackrel{\delta}{\lra}& \Lambda^2 {\cal L}_1(F)
\end{array}
$$
For a number field $F$ we have $K_2(F)_\Q=0$ Thus
 the horizontal arrows are epimorphic $\otimes \Q$. Further, for a number field 
$K_3^{\rm ind}(F)_\Q = K_3(F)_\Q$, and so the kernels of both horizontal 
arrows are isomorphic to $K_3(F)_\Q$. Thus the left arrow is an isomorphism. The lemma is proved. 
\vskip 3mm

\begin{lemma} \label{10.15.08.10} There is a canonical isomorphism
\begin{equation} \label{10.15.05.9}
{\cal B}_2({\cal O}_{F,S})_\Q \stackrel{\sim}{\lra} 
{\cal L}_{2}({\cal O}_{F,S}).
\end{equation}
\end{lemma}

{\bf Proof}. By the very definition the category of mixed Tate motives over 
${\rm Spec}{\cal O}_{F,S}$ is a subcategory of the one over ${\rm Spec}F$. Thus there is 
an inclusion of the Lie coalgebras ${\cal L}_{\bullet}({\cal O}_{F,S}) \hra 
{\cal L}_{\bullet}(F)$. It provides a commutative diagram
$$
\begin{array}{ccccccccc}
0 & \lra & K_3({\cal O}_{F,S})_\Q & \lra & {\cal L}_{2}({\cal O}_{F,S}) & \lra & \Lambda^2 {\cal L}_{1}({\cal O}_{F,S}) & \lra & 0\\
& & \downarrow & & \downarrow  & & \downarrow  & & \\
0 & \lra & K_3(F)_\Q & \lra & {\cal L}_{2}(F) & \lra & \Lambda^2 {\cal L}_{1}(F) & \lra & 0
\end{array}
$$
where the horizontal lines are exact since the fundamental Lie algebras are free, and their 
first cohomology  in the degree $2$ are given by $K_3$ of the corresponding scheme. 
Since $K_3(F)_\Q = K_3({\cal O}_{F,S})_\Q$, 
${\cal L}_{2}({\cal O}_{F,S})$ is the biggest subspace of ${\cal L}_{2}(F)$ 
which maps by the cobracket to $\Lambda^2 {\cal L}_{1}({\cal O}_{F,S})$. The lemma follows.

\subsection{A reciprocity law for $K_3^M$ of  
the function field of an elliptic curve (\cite{G1})} 
{\it The Chow dilogarithm}. Let $f_1,f_2,f_3$ be arbitrary 
rational functions on a complex curve $X$. 
Set 
$$
r_2(f_1,f_2,f_3):= 
$$
$$
{\rm Alt}_3 \Bigl(\frac{1}{6}\log|f_1| d \log|f_2| 
\wedge d\log|f_3|-\frac{1}{2}
\log|f_1| d\arg f_2 \wedge d\arg f_3 \Bigr) 
$$
 where 
${\rm Alt}_3$ is the alternation of $f_{1}, f_{2}, f_{3}$. The crucial property of this form is the following: 
$$
d r_2(f_1,f_2,f_3) = {\rm Re}\Bigl(d \log f_1 \wedge d \log f_2 \wedge d \log f_3\Bigr).
$$
Then   the  integral 
\begin{equation} \label{chowh1}
\frac{1}{2\pi i}\int_{X(\C)}r_2(f_1,f_2,f_3)
\end{equation}
converges, and provides a homomorphism 
\begin{equation} \label{chowh}
\Lambda^3 \C(X)^* \to \R,\quad  f_1\wedge f_2\wedge f_3 \lms 
\mbox{the integral (\ref{chowh1})}.
\end{equation}

Let $X$ be a regular  curve over an algebraically closed field
$k$ and $F:= k(X)^{\ast}$.  Let us define a morphism of complexes
$$
\begin{array}{ccc}
{\cal B}_2(F )\otimes
F^{\ast}&
\stackrel{\delta_2 \wedge {\rm Id}}{\longrightarrow}&\Lambda^3F^{\ast}\\
  \downarrow {\rm Res}&&\downarrow {\rm Res}\\
{\cal B}_2(k )&\stackrel{\delta_2}{\longrightarrow}& \Lambda^2k^{\ast}
\end{array}
$$
Here ${\rm Res}:= \sum_x {\rm res}_x$, where ${\rm res}_x$ is
the local residue homomorphism for the valuation $v_x$ on $F$ corresponding to a
point $x$ of $X(k)$. The local residue maps are defined as follows. 
We have ${\rm res}_x(\{y\}_2 \otimes z)=0$ unless $v_x(y) =0$. 
In the latter case 
${\rm res}_x(\{y\}_2 \otimes z) = v_x(z) \{\overline y\}_2$, 
where $\overline y$ 
denotes  projection of $y$ to the residue field of $F$ 
for the valuation $v_x$. 
The map ${\rm res}_x: \Lambda^3F^* \lra \Lambda^2k^*$ is determined uniquely by the following two conditions: 
${\rm res}_x(z_1\wedge z_2\wedge z_3)=0$ 
if $v_x(z_1) = v_x(z_2)= v_x(z_3)=0$ and 
${\rm res}_x(z_1\wedge z_2\wedge z_3)= \overline z_2 \wedge \overline z_3$ 
if $v_x(z_1) = 1$, $v_x(z_2)= v_x(z_3)=0$.

Let us present the projective plane 
${\Bbb P}^2$ as the projectivisation of the vector space $V_3$. 
A linear functional $l \in V_3^*$ is called a linear 
homogeneous function on ${\Bbb P}^2$. 
If $X$ is an elliptic curve in  
${\Bbb P}^2$, 
any rational function $f$ on $X$ can be written as a ratio
of products of {\it linear} homogeneous functions: 
$$
f = \frac{l_1\cdot ... \cdot l_k}{l_{k+1}\cdot ... \cdot l_{2k}}.
$$
This is checked by induction on the total number of 
zeros and poles using that $X$ is degree three plane curve. 

Let us return to an arbitrary curve $X$. 
For three points $a,b,c$ and a 
divisor $D =\sum n_i (x_i)$ on a 
line set 
$$
\{r(a,b,c,D)\}_2:= \sum_i n_i \{r(a,b,c,x_i)\}_2.
$$ 
Here $r(a,b,c,d)$ is the cross-ratio of four points on $P^1$, 
normalized so that $r(\infty, 0, 1, x) =x$. 

 Let  $L_i$ be  the line $l_i =0$ in the plane, 
$D_i$ the divisor $L_i \cap X$, and 
 $l_{ij}:= L_i \cap L_j$. 
The following is Theorem 6.14 in \cite{G1}. 
\begin{theorem} \label{reh}
Let $k$ be an algebraically closed field. 
Let $X$ be either an elliptic curve, or a rational curve over 
$k$, and $F:= k(X)^{\ast}$.
Then there is a  group  homomorphism 
$
h:\Lambda^3F^* \rightarrow
{\cal B}_2(k)
$ 
satisfying the following conditions: 

a) $h(k^* \wedge \Lambda^2F^*) =0$, and the following diagram is commutative:
\begin{equation} \label{6.17.02.3}
\begin{array}{ccc} \label{reh1}
{\cal B}_2(F)\otimes F^{\ast} 
& \stackrel{\delta_3}{\longrightarrow} & \Lambda^3F^{\ast}\\
&&\\
{\rm Res}\downarrow&h \swarrow&\downarrow {\rm Res}\\
&&\\
{\cal B}_2(k)&\stackrel{\delta_2}{\longrightarrow}&\Lambda^2 k^{\ast} 
\end{array}
\end{equation}

b) If $X$ is an elliptic curve,  
then for any linear homogeneous 
functions $l_0,...,l_3$ one has  
\begin{equation} \label{hrule}
h(l_1/l_0 \wedge l_2/l_0 \wedge l_3/l_0) = -\sum_{i=0}^3(-1)^i\{r(l_{i0},..., \widehat l_{ii}, ..., \l_{i3}, D_i)\}_2.
\end{equation}

c) If $X = {\Bbb P}^1$,  $t$ is a natural parameter on it, and 
$a_i \in {\Bbb P}^1(k)$, then 
\begin{equation} \label{hruleq}
h(\frac{t-a_1}{t-a_0} \wedge \frac{t-a_2}{t-a_0} \wedge \frac{t-a_3}{t-a_0} ) = -\{r(a_0, a_1, a_2, a_3)\}_2.
\end{equation}

d) If $X$ is a defined  over $\C$
then 
\begin{equation} \label{homotq}
\frac{1}{2\pi i}\int_{X(\C)}r_2(f_1\wedge f_2 \wedge f_3) =   {\cal L}_2 \Bigl(h(f_1\wedge
f_2 \wedge f_3)\Bigr).
\end{equation}

e) The map $h$ is  ${\rm Gal}(F/k)$-invariant.
\end{theorem}

\vskip 3mm

{\bf Remark}. By  Suslin's 
 reciprocity law for $K^M_3(F)$, the projection of 
${\rm Res}(\Lambda^3F^{\ast}) \subset \Lambda^2k^*$ to $K_2(k)$ is zero. 
Since by Matsumoto's 
theorem  $K_2(k)={\rm Coker}(\delta_2)$,  one has 
${\rm Res}(\Lambda^3F^{\ast}) \subset {\rm Im}(\delta_2)$. However  
${\rm Ker}(\delta_2)$ is nontrivial, so  the problem is to
 lift naturally  the map ${\rm Res}$  to a map $h$.
 
\vskip 3mm
Let $P_E$ be the abelian group of principal divisors on $E$. 
 Theorem \ref{reh} provides us a map
$$
h: \Lambda^3P_E \lra {\cal B}_2(k).
$$

\begin{lemma} \label{6.11.05.123} Let $E$ be an elliptic curve over 
an arbitrary field $k$. 
Let $D_1, D_2, D_3$ be principal divisors rational over an extension $k'$ 
of $k$. Then 
$h(D_1 \wedge D_2 \wedge D_3) \subset {B}_2(k')_\Q$. 
\end{lemma}

{\bf Proof}. Let $D= \sum_i n_i(x_i)$ 
be a divisor  on $E$ as in the Lemma. 
We can decompose it into a fraction of products 
of linear homogeneous functions $l_i$ defined over $k'$ as follows. 
Let $l_{x,y}$ (resp. $l_{x }$)  be a linear homogeneous equation of the line in ${\Bbb P}^2$ through the points $x$ and $y$ on $E$  (resp. $x $ and $-x $).  The divisor of the function $l_{x,y}/l_{x+y}$ is $(x) +(y) - (x+y) - (0)$. 
If  $D = (x)+(y) + D_1$, we write $f =
l_{x,y}/l_{x+y}\cdot f'$, so $(f') = (0) + (x+y) + D_1$. 
 After a finite number of such steps we get the desired decomposition. 
It follows that in formula (\ref{hrule}) the cross-ratios in the right hand side lie in the field $k'$. The lemma is proved.

\subsection{The modular and elliptic units revisited} 
Let $E$ be an elliptic curve over an arbitrary field $k$, and $J(k)$  the group of $k$-points of the Jacobian $J$ of $E$. 
We defined in Section 2 of \cite{GL} an abelian group $B_2(E)(k)$, usually 
denoted simply by 
$B_2(E)$, which 
has the following properties:

a) There is an exact sequence of abelian groups
\begin{equation} \label{11.3.04.1s}
0 \lra k^* \lra B_2(E) \stackrel{p}{\lra} S^2J(k) \lra 0.
\end{equation}

b) There is a canonical (up to a choice of 
a sixth root of unity) surjective 
homomorphism
$$
h: \Z[E(k)] \lra B_2(E)
$$
whose projection to $S^2J(k)$ is given by $\{a\} \lms a\cdot a$. 

c) The group $B_2(E)$ is functorial: 
an automorphism $A$ of $E$ induces an automorphism of the extension 
(\ref{11.3.04.1s}),  
which acts as the identity on the subgroup $k^*$, and whose action on $S^2J(k)$ is induced by  the map 
$A: J(k) \lra J(k)$. 

d) If $k$ is a local field then there is a homomorphism
$
H_k: B_2(E) \lra \R
$ 
whose restriction to the subgroup $k^*$ is given by $a \lms - \log|a|$.

The group $B_2(E)$ is a motivic avatar of the theta-functions 
on $E$.

\begin{corollary} \label{11.3.04.10} Let $a$ be
 an $N$-torsion element of $E(k)$. Then there is a well defined element 
$$
\theta_E(a) \in \overline k^* \otimes \Z[\frac{1}{N}].
$$

For any automorphism $A$ of $E$ we have $\theta_E(A(a)) = \theta_E(a)$.   
 \end{corollary}

{\bf Proof}. The element $p\circ h(a)$ is annihilated by $N$. 
Thus we have
$
Nh(a) \in k^*
$. 
Taking its $N$-th root we get $\theta_E(a)$. The second claim follows immediately from the property c). The corollary is proved.

\vskip 3mm
We call $\theta_E(a)$ the elliptic 
unit corresponding to the torsion point $a$ of $E$. If $E$ is defined over $\C$, it is given by 
$\theta_q(z)$ in (\ref{eq1}), where $q$ is the modulus of $E$, and $z$ is an 
$N$-torsion point on $E$. 

The following result for the first statement 
follows from Theorem 4.1 (or Corollary 4.3)
 in \cite{GL}; for the second statement see   
Theorem 4.1 on page 43 of  \cite{KL}.

\begin{lemma} \label{11.7.04.10q} 
The elliptic units $\theta_E(a)$ satisfy the following distribution relations: 
Given an isogeny $\psi: E\to E'$, one has 
\begin{equation} \label{11.8.04.4} 
\prod_{\psi(t)=t'}\theta_E(t) = \theta_{E'}(t'), \quad t' \not = 0; \qquad 
\end{equation}
\begin{equation} \label{11.8.04.5}
\prod_{\psi(t)=0}\theta_E(t) = 
\Bigl(\frac{\Delta_{E'}}{\Delta_E}\Bigr)^{1/12}.
\end{equation}
\end{lemma}
The $12$-th root of the $\Delta$-function $\Delta(\tau)$ is 
defined canonically as
$
\Delta(\tau)^{1/12} = 2\pi i ~e^{2\pi i \tau/12}\prod_{n=1}^\infty(1-q^n)^2.
$


\subsection{The  elements $\theta_E(a_1, a_2, a_3)$} Let $E$ be an 
elliptic curve over an algebraically closed field
$k$. 
For any $N$-torsion points $a,b$ on $E$ the divisor $N(\{a\} - \{b\})$ is principal, so there exists a  (non zero) 
function $f_{a,b}$ on $E$ such that ${\rm div}f_{a,b} = N(\{a\} - \{b\})$. 
It is well defined up to a non zero constant factor. 

\begin{definition} \label{6.14.05.1abv}
Let $a_1, a_2, a_3$ be $N$-torsion points on $E$. 
Then 
\begin{equation} \label{6.14.05.1a}
\theta_E(a_1: a_2: a_3):= 
-\frac{1}{N^5}\sum_{x\in E[N]}h(f_{a_1,x}, f_{a_2,x}, 
f_{a_3,x}) 
\in {\cal B}_{2}(k)_\Q
\end{equation}
where we sum over all $N$-torsion points of $E$. 
\end{definition}
Clearly this element is invariant under 
the shift $a_i \lms a_i+a$ and skew symmetric in $a_i$. 
We will show in Lemma \ref{6.28.06.2} 
that this definition does not depend on the choice of $N$.

We extend the definition of the elements $\theta_E(a_1: a_2: a_3)$ 
by linearity to a map 
$$
\theta_E: \Lambda^3\Z[E[N]] \lra {\cal B}_{2}(k)_\Q.
$$
The restriction of this map to degree zero $N$-torsion 
divisors on $E$ is given by a simpler formula: 
\begin{lemma} \label{formula1} One has 
$$
\theta_E(\{a_1\}-\{b_1\}: \{a_2\}-\{b_2\}: \{a_3\}-\{b_3\} = 
-\frac{1}{N^3} h(f_{a_1, b_1}, f_{a_2, b_2}, 
f_{a_3, b_3}). 
$$
\end{lemma}

{\bf Proof}. Clearly, $f_{a_1, b_1} = f_{a_1, x} - f_{b_1, x}$.  
The lemma follows immediately from this. 
\vskip 3mm

Observe that given a triple of elements $a_1, a_2, a_3$ 
of an abelian group $A$ 
such that $a_1+a_2 +a_3 =0$, there are elements $b_1, b_2, b_3 \in A$ such that $a_i = b_{i+1}-b_i$, where the index $i$ is modulo $3$. The triple $(b_1, b_2, b_3)$ is defined uniquely up to a shift $b_i \lms b_i +b$. 
So if $(a_1, a_2, a_3)$ are torsion points of $E$ such that $a_1+a_2 +a_3 =0$, 
we can define an element  $\theta_E(a_1,a_2,a_3)$ such that 
$$
\theta_E(b_2-b_1, b_3-b_2, b_1-b_3) = \theta_E(b_1:b_2:b_3).
$$

{\bf Remark}. 
Below we add to the elements $\theta_E(a)$, defined for $a \not = 0$, a 
formal variable $\theta_E(0)$. Therefore, strictly speaking, 
in the formulas below 
we work in $\Lambda^2(k^* \otimes \Q \oplus \Q)$ where the 
element $1$ in the extra summand $\Q$ 
corresponds to $\theta_E(0)$.
However one easily checks that in the final formulas there will be no $\theta_E(0)$. 

\vskip 3mm
The basic properties of the elements $\theta_E(a_1, a_2, a_3)$ 
are listed in the theorem below. 
\begin{theorem} \label{21:56} 
a) The differential of the element $\theta_E(a_1, a_2, a_3)$ in the
 Bloch complex 
is given by 
\begin{equation} \label{22:10}
\delta_2: \theta_E(a_1, a_2, a_3) \lra \theta_E(a_1)\wedge \theta_E(a_2) 
+ \theta_E(a_2)\wedge \theta_E(a_3) + \theta_E(a_3)\wedge \theta_E(a_1).
\end{equation}

b) The elements $\theta_E(a_1, a_2, a_3)$ satisfy
 the following relations: 

\begin{itemize}
\item
The dihedral symmetry relations:
\begin{equation} \label{22:101}
 \theta_E(a_1, a_2, a_3) = (-1)^
{{\rm sign}(\sigma)}\theta_E(a_{\sigma(1)}, 
a_{\sigma(2)}, a_{\sigma(3)}), \qquad \sigma \in S_3, 
\end{equation}
\begin{equation} \label{22:102}
 \theta_E(a_1, a_2, a_3) = \theta_E(\varepsilon (a_{1}), \varepsilon (a_{2}), \varepsilon (a_{3})), \qquad \varepsilon \in {\rm Aut}(E).
\end{equation}

\item
The distribution relations: 
Given an isogeny $\psi: E\to E'$, one has 
\begin{equation} \label{11.8.04.4a} 
\sum_{\psi(a_i)=a_i'}\theta_E(a_1, a_2, a_3) = 
\theta_{E'}(a_1', a_2', a_3'), \quad a_i' \not = 0; \qquad 
\end{equation}
\end{itemize}
\end{theorem}

{\bf Proof}.   a) The formula is equivalent to the following one, 
 which we are going to prove:
\begin{equation} \label{22:10j}
\delta_2: \theta_E(a_1: a_2: a_3) \lra 
\end{equation}
$$
\theta_E(a_1-a_2)\wedge \theta_E(a_2-a_3) 
+ \theta_E(a_2-a_3)\wedge \theta_E(a_3-a_1) + \theta_E(a_3-a_1)\wedge \theta_E(a_1-a_2).
$$
Denote the right hand side of (\ref{22:10j}) by $C(a_1:a_2:a_3)$. 
Let us prove first that 
\begin{equation} \label{22:10jj}
\delta_2 \circ \theta_E(\{a_1\} - \{x\}:\{a_2\} - \{x\}: \{a_3\} - \{x\}) 
= 
\end{equation}
$$
C(a_1:a_2:a_3) - C(x:a_2:a_3) +  C(a_1:x:a_3) 
- C(a_1:a_2:x). 
$$
Using  
$f_{a_i,x} = \theta_E(t-a_i)^N/\theta_E(t-x)^N$, we get 
$$
\frac{1}{N^3}{\rm Res}
\left( f_{a_1,x} \wedge f_{a_2,x} \wedge f_{a_3,x}\right) = C(a_1:a_2:a_3) - C(x:a_2:a_3) +  C(a_1:x:a_3) 
- C(a_1:a_2:x). 
$$
Thanks to (\ref{6.17.02.3}), this implies (\ref{22:10jj}). 
Now we deduce formula (\ref{22:10j}) from this. We claim that 
\begin{equation} \label{rty}
\sum_{x\in E[N]}C(x:a_2:a_3)=0.
\end{equation} 
Having this formula, and taking the sum of 
formula (\ref{22:10jj}) over all $x \in E[N]$, we get formula (\ref{22:10j}). 

We will obtain formula 
(\ref{rty}) as a special case of the 
following simple general statement.

\vskip 3mm 
Let $A$ and $B$ be abelian groups. 
Let 
$$
\Phi: \Lambda^3 \Z[A] \lra B, \qquad \{a\} \wedge \{b\} \wedge \{c\} \lms 
\Phi(a: b: c).
$$ 
be a group homomorphism. 

\begin{lemma} \label{6.28.05.1} 
Let us assume that the group $A$ is finite and 
the map $\Phi$ satisfies the following properties: 
$$
\Phi(-a: b: c) = \Phi(a: b: c), \qquad 
\Phi(a: b: c) = \Phi(a+x: b+x: c+x) \quad \mbox{for any $x\in A$}.
$$
Then one has modulo $2$-torsion
$$
\sum_{x\in A}\Phi(a: b: x)=0.
$$ 
\end{lemma}

{\bf Proof}. Indeed, one has 
$$
\sum_{x\in A}\Phi(a: b: x) = \sum_{x\in A}\Phi(a-(a+b): b-(a+b): x-(a+b)) =
$$
$$ 
\sum_{y\in A}\Phi(-b: -a: y) = -\sum_{y\in A}\Phi(a: b:  y).
$$
The Lemma is proved.
\vskip 3mm

 Since $\theta_E(y) = \theta_E(-y)$, 
and $a_2, a_3$ are $N$-torsion points, the element $C(a_1:a_2:a_3)$ satisfies the dihedral symmetry relations, and by definition is invariant under the shift of the arguments. Thus formula 
(\ref{rty}) follows from Lemma \ref{6.28.05.1}. 
The part a) of the theorem is proved. 

b) The dihedral symmetry relations follow from Theorem \ref{reh} and 
the definition. 

Let us prove the  distribution relations. 
Thanks to the distribution relations for elliptic units, the difference between the left and the 
right hand sides of (\ref{11.8.04.4a}) is killed by the differential 
$\delta_2$. Since the kernel of the differential $\delta_2$ on the group ${\cal B}_2$ is rigid by the very definition, the difference is 
constant on the modular curve: Indeed, take the  modular curve as the curve 
 $X$ in the definition in Section 2.1. So it remains to check that it is zero at a single point. As such one can take a CM point corresponding to an elliptic curve 
with the endomorphism ring ${\cal O}_K$. 
Then the difference lies in 
${\rm Ker}\delta_2 \subset {\cal B}_2(K')$ 
for a certain number field $K'$ extending $K$ (use Lemma \ref{6.11.05.123}). 
Thanks to the 
injectivity of the regulator map, it is enough to show that the 
dilogarithm kills this element and its Galois conjugated. 
Now the claim follows from (\ref{homotq}) plus the fact that 
the Chow dilogarithm satisfies the distribution relations on the nose. The claim is proved. The theorem is proved. 

\vskip 3mm
{\bf Remark}. Observe that 
$\theta_E(a, -a, 0) =0$ since by the dihedral relations 
one has $\theta_E(a, -a, 0) = - \theta_E(-a, a, 0) = -\theta_E(a,  -a, 0)$. 
Formula (\ref{22:10}) makes sense if $a_1=0$. Indeed, thanks to $\theta_E(a)= \theta_E(-a)$, 
it gives 
$$
\delta_2 \circ \theta_E(0, a, -a) = \theta_E(0) \wedge \theta_E(a) + 
\theta_E(a) \wedge \theta_E(-a) + \theta_E(-a) \wedge \theta_E(0) =  0.
$$

\begin{lemma} \label{6.28.06.2}
The element (\ref{6.14.05.1a}) does not depend on the choice of $N$. Namely,  
if $N|M$, then the elements defined by formula (\ref{6.14.05.1a}) 
using $N$ or $M$ coincide. 
\end{lemma}

{\bf Proof}. For $M$-torsion points $a,b$, denote by 
$f^{(M)}_{a, b}$ a function with the divisor $M(\{a\}-\{b\})$. 
Given $M$-torsion points $a_1, a_2, a_3, x$
one has, by Lemma \ref{formula1},
$$
-\frac{1}{M^3} h(f^{(M)}_{a_1, x}, f^{(M)}_{a_2, x}, 
f^{(M)}_{a_3, x})  = \theta_E(a_1: a_2: a_3) - \theta_E(x: a_2: a_3) + 
\theta_E(a_1: x: a_3) - \theta_E(a_1: a_2: x). 
$$
Assume that $a_i, x$ are $N$-torsion points. 
Averaging over $x \in E[N]$, and using Lemma 
\ref{6.28.05.1},  we get 
$$
\frac{1}{N^2}\sum_{x \in E[N]} -\frac{1}{M^3} h(f^{(M)}_{a_1, x}, f^{(M)}_{a_2, x}, 
f^{(M)}_{a_3, x})  = \theta_E(a_1: a_2: a_3). 
$$
Since $(f^{(M)}_{a_i, x})^{M/N} = f^{(N)}_{a_i, x}$, the lemma follows. 
\vskip 3mm

Below we will need the following 
 special case of the definition of elements 
$\theta_E(a:b:c)$ for the modular curve $Y_1(N)$. 
Recall that points of $Y_1(N)$ parametrize pairs $(E, p)$, where 
$p$ is an $N$-torsion point of $E$, generating the subgroup $E[N]$. 
So for any residues $\alpha_i \in \Z/N\Z$ there are elements 
$\theta_E(\alpha_1 p: \alpha_2 p: \alpha_3 p)$. 

\begin{lemma} \label{6.28.06.3} One has 
$$
\theta_E(\alpha_1 p: \alpha_2 p: \alpha_3 p) = 
\frac{1}{N}\sum_{\beta \in \Z/N\Z} 
h(f^{(N)}_{\alpha_1 p, \beta p}, f^{(N)}_{\alpha_2 p, \beta p}, 
f^{(N)}_{\alpha_3 p, \beta p}).
$$
\end{lemma}

{\bf Proof.} Follows from Lemmas \ref{6.28.05.1} 
and  \ref{formula1} the same way as 
Lemma \ref{6.28.06.2}. 

\vskip 3mm

\subsection{Euler complexes on modular curves} 
Recall the full level $N$  modular curve $Y(N)$. It 
 is defined over the cyclotomic field $\Q(\mu_N)$. 
Denote by $F_{Y(N)}$ its function field. 
Let ${\cal E}_N$ be the universal elliptic curve 
over  $Y(N)$. 

Recall that for a torsion section $a$, the element 
 $\theta_E(a)$ is a modular unit, i.e. lies in  
 ${\cal O}^*_{Y(N)}\otimes \Z[\frac{1}{N}]$. 
\begin{definition} \label{6.14.05.3}
The  $\Q$-vector space ${\Bbb E}^1_N$ is 
 generated by 
the elements 
$$
\theta_E(a_0:a_1:a_2) \in {\cal B}_2(F_{Y(N)})_\Q
$$ when $a_i$'s run 
through all $N$-torsion points of the universal elliptic curve 
${\cal E}_N$.

The $\Q$-vector space ${\Bbb E}^2_N$ is generated by the elements  
$$
\theta_E(a) \wedge\theta_E(b)\in 
\Lambda^2{\cal O}^*_{Y(N)}\otimes \Q
$$
when $a,b$ run through all non-zero $N$-torsion points of 
the universal elliptic curve ${\cal E}_N$. 
\end{definition}

It follows from the formula (\ref{22:10jj}) that we get a complex 
${\Bbb E}^{\bullet}_N$, placed in degrees $[1,2]$:
$$
{\Bbb E}^{\bullet}_N: \quad {\Bbb E}^1_N\stackrel{\delta_2}{\lra} {\Bbb E}^2_N.
$$ 
We  call it 
the {\it Euler complex} of the modular curve $Y(N)$. 
\begin{lemma} \label{formula2}
For any three $N$-torsion points $a,b,c$ of $E$ one has 
$$\theta_E(a: b: c) \in {\cal B}_2(Y(N))\otimes \Q.$$ 
\end{lemma}

{\bf Proof}. Follows from the definition of ${\cal B}_2(Y(N))$ 
formula (\ref{22:10jj}) and the fact that 
$\theta_E(a)\in {\cal O}^*_{Y(N)}$.
\vskip 3mm

Thus the Euler complex is a subcomplex of the Bloch complex 
of the modular curve $Y(N)$:
$$
\begin{array}{ccc}
{\Bbb E}^1_N& \lra &{\Bbb E}^2_N\\
\downarrow &&\downarrow \\
\Bigl({\cal B}_2(Y(N))&\lra& \Lambda^2{\cal O}^*_{Y(N)}\Bigr)_\Q.
\end{array} 
$$
Since ${\cal O}^*_{Y(N)}\otimes \Q$ is generated by the 
elements $\theta_E(a)$, when $a \in E[N]$, the right arrow in the diagram above is an isomorphism, and the left one is inclusion by definition. 
When it is an isomorphism? In other words, when ${\cal B}_2(Y(N))\otimes \Q$ is 
generated over $\Q$ by the elements $\theta_E(a, b, c)$, where $a,b,c \in E[N]$?

Recall the modular complex 
\begin{equation} \label{4-12.1z}
{\rm M}_{(2)}^{\ast}:= {\rm M}_{(2)}^{1} \lra {\rm M}_{(2)}^{2}.
 \end{equation} 
Let $\Gamma$ be a subgroup of $GL_2(\Z)$. Projecting the modular 
complex onto the modular curve $Y_{\Gamma}:=  {\cal H}_2/\Gamma$, we get the modular triangulation of the latter. 
Its chain complex is identified with ${\rm M}_{(2)}^{\ast}\otimes_\Gamma\Q$. 

\begin{proposition} \label{6.16.05.1} 
For every positive integer 
$N$ there is a canonical surjective homomorphism of complexes
\begin{equation} \label{6.16.05.s1} 
\begin{array}{ccc} 
{\rm M}_{(2)}^{1}&\lra&{\rm M}_{(2)}^{2}\\
\lambda_N^1\downarrow & &\downarrow \lambda_N^2\\
{\Bbb E}^{1}_N&\lra&{\Bbb E}^{2}_N
\end{array}
\end{equation}
providing a map of the modular complex for $Y(N)$ to the corresponding Euler complex:
\begin{equation} \label{6.16.05.po} 
\lambda_N^*: {\rm M}_{(2)}^{\ast}\otimes_{\Gamma(N)} \Q\lra {\Bbb E}^{*}_N.
\end{equation}
\end{proposition} {\bf Proof}. 
The modular complex for $Y(N)$ has the following description: 
$$
\Z[\Gamma(N)\backslash GL_2(\Z)] \otimes _{D^1}\chi_1 \lra 
\Z[\Gamma(N)\backslash GL_2(\Z)] \otimes _{D^2}\chi_2. 
$$
Here $D^1$ is the subgroup of order $12$ 
of $GL_2(\Z)$ stabilizing the modular triangle 
with vertexes $(0,1,\infty)$. 
Further, $D^2$ is the subgroup of order $8$ of $GL_2(\Z)$ 
stabilizing the geodesic 
 $(0,\infty)$.  The characters $\chi_1$ and 
$\chi_2$ are given by the determinant. 
Observe that $\Gamma(N)\backslash GL_2(\Z)= GL_2(\Z/N\Z)$. 
A point of the universal elliptic curve ${\cal E}_N$ is given by a triple 
$(E; p_1, p_2)$, where $E$ is an elliptic curve, and $(p_1, p_2)$ is a 
basis in the abelian group $E[N]$. The maps $\lambda_N^1$ and $\lambda_N^2$ act 
on the generator parametrized by a matrix 
$$
\left 
(\matrix{a&b \cr  c& d\cr}\right ) \in GL_2(\Z/N\Z)
$$
as follows: 
\begin{equation} \label{6.16.05.s2}
\lambda_N^1: \left(\matrix{a&b \cr  c& d\cr}\right ) \lms 
\theta_E(ap_1 + cp_2, bp_1+ dp_2, 
-(a+b)p_1 -(c+d)p_2); \qquad  
\end{equation}
\begin{equation} \label{6.16.05.s3}
\lambda_N^2: \left 
(\matrix{a&b \cr  c& d\cr}\right ) \lms \theta_E(ap_1 + cp_2) \wedge 
\theta_E(bp_1 + dp_2).
\end{equation}
Now the claim that the map (\ref{6.16.05.s2}) gives rise to a group homomorphism 
$\lambda_N^1$ given by the left 
arrow in (\ref{6.16.05.s1}) is equivalent to the dihedral symmetry relations 
(\ref{22:101}). 
The claim that the map (\ref{6.16.05.s3}) gives rise to a group 
homomorphism $\lambda_N^2$
given by the right  
arrow in (\ref{6.16.05.s1}) follows from $\theta_E(-a) = 
\theta_E(a)$ and the skew-symmetry of the wedge product.  
The proposition is proved. 
\vskip 3mm
Summarizing, we arrive at the diagram on Fig \ref{ec1}. 

\begin{corollary} The map $\lambda_N^2$ gives rise to a canonical homomorphism
$$
\lambda_N: H^1(\Gamma(N), \Z) \lra K_2(Y(N))\otimes \Z[\frac{1}{N}].
$$
\end{corollary}

\vskip 3mm
The $H^2$ cohomology groups of the Euler complexes
for various $N$ deliver the 
Beilinson-Kato Euler system on the tower of modular curves. 
Let us explain this in more detail. 

The {\it Euler complex datum} 
on a modular curve $Y(N)$ is a homomorphism of complexes (\ref{6.16.05.po}). The Euler complex is the image of this homomorphism. The source 
of the map (\ref{6.16.05.po}) does not depend on the 
choice of the modular curve. Therefore each element $\gamma$ 
of the modular complex gives rise to a collection of elements 
$\lambda_N(\gamma)$ in Euler complexes. In particular, if 
$\gamma$ is the element corresponding to the geodesic $(0, \infty)$ on 
${\cal H}_2$, the elements $\lambda_N(\gamma)$, projected to 
$K_2(Y(N)$, give rise to 
an Euler system in $K_2$ --  this 
is deduced from the results of Kato \cite{Ka}.

\section{Specialization of Euler complexes 
at a cusp}

\subsection{Some identities in the Bloch group}
Let $k$ be an arbitrary field. We define the motivic dilogarithm 
$
{\rm Li}^{\cal M}_{2}(x) := \{x\}_2 \in { B}_2(k)$. 
Following Section 2.1 of \cite{G2}, 
the motivic double logarithm is the following element of the 
Bloch group:
$$
{\rm Li}^{\cal M}_{1,1}(x,y):= \{\frac{xy-y}{1-y}\}_2 - 
\{\frac{y}{y-1}\}_2 - \{xy\}_2 \in B_2(k).
$$
This definition was suggested by a similar
 identity between the corresponding multivalued analytic functions, 
easily proved by differentiation, see loc. cit..  

It is easy to check the crucial formula (\ref{6.14.05.5qaz}) 
for the coproduct $\delta_2{\rm Li}_{1,1}(a,b)$, where $a^N = b^N =1$. 

Let us introduce a standard homogeneous notation for the double logarithm: 
\begin{equation} \label{formula3}
{\rm I}^{\cal M}_{1,1}(a_1: a_2: a_3) := 
{\rm Li}^{\cal M}_{1,1}(\frac{a_2}{a_1}, \frac{a_3}{a_2}).
\end{equation}
It is a motivic version of the iterated integral $I(a_0:a_1:a_2) = \int_0^{a_2}
\frac{dt}{t-a_1}\circ \frac{dt}{t-a_2}$. 
\begin{lemma} \label{6.26.05.1} Let $a_i^N=1$. 
Then one has the following identity in 
 $B_2(k)$ modulo $6$-torsion:
$$
\sum_{i=0}^3(-1)^i{\rm I}^{\cal M}_{1,1}(a_0: \ldots :\widehat a_i: \ldots  : a_3) = \{r(a_0, a_1, a_2, a_3)\}_2.
$$
\end{lemma}

{\bf Proof}. One has 
$$
{\rm I}^{\cal M}_{1,1}(a_1: a_2: a_3) = \{\frac{a_3(a_2-a_1)}{a_1(a_2-a_3)}\}_2  - 
\{\frac{a_3}{a_3-a_2}\}_2 - \{\frac{a_3}{a_1}\}_2 
$$
After the alternation the last two terms on the right in this formula will
 disappear. Since 
$$
r(0, a_2, a_3, a_1) = \frac{a_3(a_2-a_1)}
{a_1(a_2-a_3)}
$$
 and $\{r(a_{\sigma(1)}, ..., a_{\sigma(4)}\}_2 = (-1)^{|\sigma|}
\{r(a_{1}, ..., a_{4}\}_2$ modulo $6$-torsion, 
what remains is the five-term relation. 
The lemma is proved. 

\begin{lemma} \label{8.29.06.1}
Let $a^N_i=1$. 
Then one has  in $B_2(k)\otimes \Q$ 
\begin{equation} \label{8.29.06.2}
\frac{1}{N}\sum_{x^N=1} \{r(a_1, a_2, a_3, x)\}_2 = 
{\rm I}^{\cal M}_{1,1}(a_1 : a_2:  a_3).
\end{equation}
\end{lemma}

{\bf Proof}. Lemma \ref{6.28.05.1} implies that 
$$
\sum_{x^N=1} {\rm I}^{\cal M}_{1,1}(a_1 : a_2:  x) = 0.
$$
This and Lemma \ref{6.26.05.1} imply the identity. The Lemma is proved.

\subsection{The specialization  map} 
Let ${\cal O}$ be a discrete valuatiuon ring with a uniformizer $q$, 
the residue field $k$, the valuation homomorphism $v$, and 
the fraction field $F$. 
There is a canonical homomorphism ${\cal O} \to k, f \lms \overline f$. 
We define a specialization map 
$$
{\rm Sp}_q: F^* \lra k^*, \quad 
f \lms \overline {f/q^{v(f)}}.
$$

\begin{lemma} \label{Commute} The  specialization maps 
$$
{\rm Sp}_q\{f\}_2 = \left\{ \begin{array}{ll}
\overline f
& \quad \mbox{if } f \in {\cal O} ^*\\ 
 0 &  
\quad \mbox{otherwise}   \end{array} \right., \qquad 
{\rm Sp}_q(f_1 \wedge   f_2):= 
{\rm Sp}_qf_1 
\wedge {\rm Sp}_q f_2
$$ 
gives rise to a  homomorphism of complexes, modulo $2$-torsion,
\begin{equation}
\begin{array}{ccc}
{\cal B}_2(F)& \lra &\Lambda^2F^*\\
\downarrow {\rm Sp}_q&& \downarrow {\rm Sp}_q\\
{\cal B}_2(k)&\lra & \Lambda^2k^*
\end{array}
\end{equation}

\end{lemma}

{\bf Proof}. Let us show that the specialization map 
commutes with the differentials. Indeed, set $f = q^nf_0$ where 
$f_0 \in {\cal O}^*$. Then if $n \not = 0$, one has 
${\rm Sp}_q\{f\}_2 =0$. On the other hand, if $n>0$, then 
${\rm Sp}_q(1-q^nf_0)\wedge q^nf_0 = 1 \wedge \overline f_0 =0$. Similarly, 
if $n<0$, we get ${\rm Sp}_q((1-f)\wedge f) = (-\overline f_0)\wedge 
\overline f_0 =0$ modulo $2$-torsion. If $n=0$ but $\overline f_0 \not = 1$, 
the commutativity of the diagram is clear. Finally, if $n=0$ and 
$\overline f_0 = 1$, then 
$\{\overline f\}_2=0$ and ${\rm Sp}_q((1-f)\wedge f) = (1-f)/q^{v(1-f)}\wedge 1 = 0$. 
We left to the reader to check that the specialization map 
is a well defined map ${\cal B}_2(F) \lra {\cal B}_2(k)$. The lemma is proved. 
\vskip 3mm
 Therefore there  is a specialization homomorphism of complexes 
corresponding
 to a local parameter $q$ at a $k$-point  $x$ of a curve $X$ over $k$.

\subsection{Specialization of the Euler complex at a cusp on $Y_1(N)$}

Let $\infty$ be the 
cusp 
on a modular curve obtained by projection of the 
point  $\infty \in {\cal H}_2$. 
Then, for some integer $N$, 
 $q = {\rm exp}(2\pi i \tau/N)$ is a local parameter 
at this cusp. Using it, we define the specialization map  ${\rm Sp}_q$ 
 at the cusp $\infty$. 

\begin{theorem} \label{6.25.05.1}
a) The specialization map ${\rm Sp}_q$ at the cusp $\infty$ on $Y_1(N)$ 
provides a surjective homomorphism from the Euler complex on $Y_1(N)$ to the level $N$ cyclotomic complex:
\begin{equation} \label{6.25.05.2}
\begin{array}{ccc}
{\Bbb E}^1(Y_1(N))& \lra &{\Bbb E}^2(Y_1(N))\\
&&\\
\downarrow {\rm Sp}_q &&\downarrow {\rm Sp}_q\\
&&\\
{C}_2(N)&\lra &\Lambda^2{C}_1(N)
\end{array}
\end{equation}
It intertwines the maps
 of the modular complex to the Euler and cyclotomic complexes.

b) If $N=p$ is a prime number, it gives rise to an isomorphism of complexes
\begin{equation} \label{6.25.05.21}
\begin{array}{ccc}
{\Bbb E}^1(Y_1(p))& \lra &{\Bbb E}^2(Y_1(p))\\
&&\\
=\downarrow {\rm Sp}_q &&=\downarrow {\rm Sp}_q\\
&&\\
{C}_2(p)&\lra &\Lambda^2 \widehat {C}_1(p)
\end{array}
\end{equation}
\end{theorem}

{\bf Proof}. a) The specialization of  modular units 
$\theta_E(z)$ on $Y(N)$, in the notation of formula (\ref{eq1}), is:
\begin{equation} \label{6.25.05.4}
{\rm Sp}_q\theta_E(z) = \left\{ \begin{array}{ll}
(1-e^{2\pi i \alpha_2}) 
& \quad \mbox{if } \alpha_1 =0, \alpha_2\not = 0\\ 
 1 &  
\quad \mbox{otherwise}   \end{array} \right. \quad \mbox{modulo torsion}.
\end{equation}

For the modular curve $Y_1(N)$, 
there are modular units $\theta_E(\alpha p)$, where $(E, p)$ is a point of 
$Y_1(N)$, and $\alpha \in \Z/N\Z - 0$. We get
\begin{equation} \label{6.25.05.3}
{\rm Sp}_q\theta_E(\alpha p) = 
(1-e^{2\pi i \alpha})   
\quad \mbox{modulo torsion}.
\end{equation}
So the specialization of a modular unit at the cusp $\infty$ 
is a cyclotomic $N$-unit. 
This implies that the right arrow in (\ref{6.25.05.2}) is surjective.

\begin{lemma} \label{6789} One has 
$$
{\rm Sp}_q\theta_E(\alpha_1 p, \alpha_2 p,\alpha_3 p) = 
{\rm Li}^{\cal M}_{1,1}(\zeta_N^{\alpha_1}, \zeta_N^{\alpha_2}), 
\quad \zeta_N = {\rm exp}(2\pi i /N).
$$
\end{lemma}

{\bf Proof}. Both the left and the right hand sides are zero if one of $\alpha_i$ is zero modulo $N$. So we may assume $\alpha_i$  are non zero.

One can realized the family $(E, p)$ of pairs (an elliptic curve plus 
a point $p$ of order $N$ on $E$) parametrized by 
a neighborhood of a cusp $\infty$, 
as a  family of plane elliptic curves $E_q$ with an $N$-torsion 
section $p$, degenerating at $q=0$ to a nodal curve with 
the torsion point $\zeta_N$. Then $f_{\alpha_1 p, b}\wedge 
f_{\alpha_2 p, b}\wedge f_{\alpha_3 p, b}$, where $b$ is a non zero multiple of 
$p$,  degenerates, modulo $N$-torsion, 
 to an element $f_{\zeta_N^{\alpha_1}, \zeta_N^{\beta}}\wedge 
f_{\zeta_N^{\alpha_2}, \zeta_N^{\beta}}\wedge f_{\zeta_N^{\alpha_3}, \zeta_N^{\beta}}$ 
on $\Lambda^3\Q({\Bbb G}_m)^*$. It remains to apply Lemma \ref{6.28.06.3}, 
formula (\ref{hruleq}), and 
identity (\ref{8.29.06.2}). 
(The reader may compare
 this with the proof of the formula (\ref{11.8.04.4a})). The Lemma is proved. 

\vskip 3mm
It follows from Lemma \ref{6789}  that the specialization map intertwines 
the maps of the modular complex to the Euler and cyclotomic ones.

b) If $p$ is a prime, the only relations between the 
modular units, as well as the cyclotomic units,  are the parity relations 
$\theta_E(-a) = \theta_E(a)$, and similar ones for the cyclotomic case. 
Therefore we get a map of complexes which is surjective by definition, 
and is an isomorphism on the right. 
We claim that the left arrow is injective. 
The specialization map, restricted to 
the kernel of the differential in the Euler complex, is an isomorphism. 
Indeed,   
the Bloch group ${\cal B}_2(X)$ 
is rigid by its very definition. So any element of ${\rm Ker}\delta_2$ 
is constant on $X$. The specialization of 
a constant  on $X-x$ equals to its value
at any nearby point. So we get an isomorphism of complexes. 
 The theorem is proved.

\section{ Imaginary quadratic fields and 
tessellations of the hyperbolic space} 

Let ${\cal O}_K$ be the ring of integers in the imaginary quadratic field 
$K= \Q(\sqrt{-d})$. 
Following the work of Cremona \cite{C}, we 
recall the classical description (due to Bianchi \cite{Bi}) 
of the fundamental domains for the group $GL_2({\cal O}_K)$ 
acting on the hyperbolic space ${\cal H}_3$ for $d=1,2,3,7,11$.  
Using this, we describe 
explicitly the chain complex of the corresponding 
tessellation of the hyperbolic plane. 

In the next section we relate this complex with the Euler complex restricted to a CM point of the modular curve corresponding to an elliptic curve $E_K$ 
with the endomorphism ring ${\cal O}_K$. 

\subsection{Bianchi tessellations of the hyperbolic space} 
Let ${\cal H}_3$ be the three dimensional hyperbolic space. The group $GL_2(\C)$ acts 
on ${\cal H}_3$ through its quotient $PGL_2(\C)$. 
We present the hyperbolic space as 
an upper half space
$$
{\cal H}_3 = \left\{(z,t): z\in \C, t\in \R, t>0 \right \}.
$$

Let $K$ be an imaginary quadratic field $\Q(\sqrt {-d})$, 
where $d$ is a square-free positive integer. 
Let ${\cal O}_K$  be the ring of integers in $K$. 
Then $GL_2({\cal O}_K)$ acts on ${\cal H}_3$ discretely 
  through its quotient $PGL_2({\cal O}_K)$. 
Extending ${\cal H}_3$ by the set of cusps 
$K \cup \{\infty\}$ on the boundary, we form ${\cal H}_3^* := {\cal H}_3 \cup K \cup \{\infty\}$.

The kernel of the $GL_2({\cal O}_K)$-action 
on ${\cal H}_3$ is the center of $GL_2({\cal O}_K)$. It
 is isomorphic to the group of roots of unity in $K$, which coincides with 
 the units group  ${\cal O}_K^*$. 
It is the  subgroup $\{\pm 1\}$, unless $d=1, 3$, i.e. 
${\cal O}_K$ is the ring of Gaussian or Eisenstein integers, and  its order is 
$4$ or $6$ respectively. 
The group 
$SL_2({\cal O}_K)$ is the kernel of the determinant 
map $GL_2({\cal O}_K) \lra {\cal O}^*_K$. 
Observe that $PSL_2({\cal O}_K)$ is an index two 
subgroup in $PGL_2({\cal O}_K)$. 

\vskip 3mm
{\it The Bianchi tessellation and the modular complex}. 
The Bianchi tessellation 
is a classical $GL_2({\cal O}_K)$-invariant 
tessellation of ${\cal H}^*_3$ on geodesic polyhedrons. 
It can be defined by the following construction, due to Voronoi. 
Let $L$ be a free rank two ${\cal O}_K$-module. Set 
$V:= L\otimes_{{\cal O}_K}\C$. 
We realize ${\cal H}_3$ as the quotient of the cone  of 
positive definite Hermitian forms on $V^*$ modulo $\R_{>0}$. 
A non-zero vector $l \in L$ gives rise to a degenerate Hermitian form 
$\varphi_l$, where $\varphi_l(f):= |\langle f, l\rangle|^2$ on $V^*$. Let us 
consider the set of forms $\varphi_l$ when $l$ run through the set of all non zero primitive vectors in the lattice $L$. 
The convex hull of this set is an infinite polyhedron. 
Its projection to ${\cal H}_3$ provides a tessellation of the latter. 
It is the Bianchi tessellation. 
The set of its vertexes coincides with the set of cusps $K \cup \{\infty\}$. 
The interiors of the cells of dimension $\geq 1$ are inside of ${\cal H}_3$. The stabilizers of the cells are 
finite.

\begin{definition} \label{11.2.04.2}
The modular complex ${\rm M}_{{\cal O}_K}^*$ for $GL_2({\cal O}_K)$ is the chain complex of the Voronoi 
tessellation of the hyperbolic space ${\cal H}_3$ for the group 
$GL_2({\cal O}_K)$:
$$
{\rm M}_{{\cal O}_K}^*: \qquad 
{\rm M}_{{\cal O}_K}^0 \stackrel{\partial}{\lra} {\rm M}_{{\cal O}_K}^1 
\stackrel{\partial}{\lra} {\rm M}_{{\cal O}_K}^2. 
$$
The extended modular complex ${\Bbb M}_{{\cal O}_K}^*$ is the chain complex of the Voronoi 
tessellation of ${\cal H}^*_3$ for the group 
$GL_2({\cal O}_K)$:
$$
{\Bbb M}_{{\cal O}_K}^*: \qquad 
{\rm M}_{{\cal O}_K}^0 \stackrel{\partial}{\lra} {\rm M}_{{\cal O}_K}^1 
\stackrel{\partial}{\lra} {\rm M}_{{\cal O}_K}^2 
\stackrel{\partial}{\lra} {\rm M}_{{\cal O}_K}^3. 
$$
Here ${\rm M}_{{\cal O}_K}^i$ 
is the group generated by the oriented $(3-i)$-cells of the 
tessellation.  
\end{definition} 
The extended modular complex is a complex of right 
$GL_2({\cal O}_K)$-modules. 

\subsection{Bianchi tessellations for Euclidean fields} 
Now let us assume that 
\begin{equation} \label{11.3.04.1}
\mbox{$K$ is one of the five Euclidean field $\Q(\sqrt {-d})$, where $d = 1, 2, 3, 7, 11.$}. 
\end{equation}
The Voronoi tessellation in this case was considered by Bianchi \cite{Bi}. 
Its main features are the following: 

\begin{itemize}
\item
 The edges of the Bianchi tessellation are obtained 
from a single geodesic ${\Bbb G}= (0, \infty)$ by the action of the group $GL_2({\cal O}_K)$. 

\item The three dimensional cells of the Bianchi tessellation 
are obtained from a single 
{\it basic geodesic polyhedron} ${\Bbb B}_d$ by the action of 
the group $GL_2({\cal O}_K)$. 
\end{itemize}

There is also a {\it basic triangle} ${\Bbb T} = (0, 1, \infty)$. The orbits of ${\Bbb T}$ give 
all the two dimensional cells of the Voronoi tessellation if and only if $d=1,3$.

Reflecting the basic polyhedron 
by the faces, and repeating this procedure 
infinitely many times 
with the obtained polyhedrons,  we recover the Bianchi tessellation of ${\cal H}_3$.
 Below we describe the basic polyhedrons for the Euclidean fields.

{\bf The Bianchi  tessellation  
 for the Gaussian integers}. 
The units group $\Z[i]^*$ is generated by $i$. 
The basic polyhedron ${\Bbb B}_1$ is the geodesic octahedron  
described by its vertexes: 
$$
{\Bbb B}_1 = (0, 1, i, 1+i, (1+i)/2, \infty).
$$

{\bf  The Bianchi  tessellation 
 for the Eisenstein integers}. 
Let $\Z[\rho]$, where 
 $\rho = \frac{1+\sqrt{-3}}{2}= {\rm exp}(\pi i/3)$, be  the ring 
of Eisenstein integers. Its units group 
is cyclic of order $6$, generated by $\rho$. 
The basic polyhedron is the geodesic tetrahedron 
$$
{\Bbb B}_3 = (0, 1, \rho, \infty).
$$

\begin{figure}[ht]
\centerline{\epsfbox{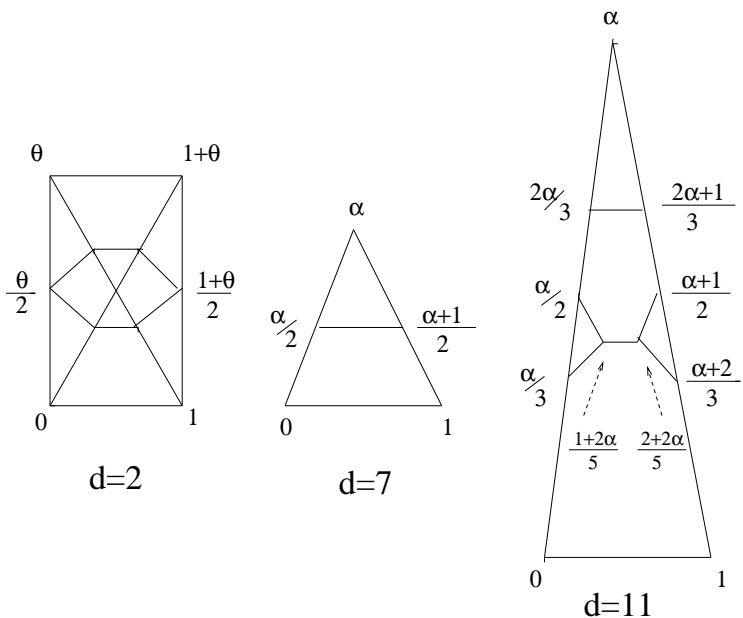}}
\caption{Plans of the basic polyhedrons for $d=2, 7, 11$}
\label{kj8}
\end{figure}

{\bf  The Bianchi  tessellations 
 for the Euclidean fields with $d=2,7,11$}. 
We present on Figure \ref{kj8} the plans of the basic polyhedrons 
${\Bbb B}_2, {\Bbb B}_7, {\Bbb B}_{11}$.  
We set there 
$$
\theta:= \sqrt{-2}, \qquad \alpha:= (1+\sqrt{-d})/2, \quad d=7,11
$$
Each of the plans shows the finite vertexes of a basic polyhedron as the vertexes of the plan, and the projections of the non-vertical 
geodesic edges of the polyhedron 
from the infinity as the edges of the plan. Observe that 
in the addition to the triangular faces we 
have quadrilateral faces for $d=2$ and $d=7$, and hexagonal faces for $d=11$. 
The action of the group ${\Bbb D}^0_d$ on the set of faces of 
${\Bbb B}_{d}$ has two orbits 
for $d=2,7,11$. The vertical faces of the polyhedrons project to the exterior sides of the plan. 
So for $d=2$ we have two triangular vertical faces, 
and two quadrilateral ones. For $d=7$ we have two quadrilateral vertical faces, 
and one triangular. For $d=11$ we have two hexagonal vertical faces, 
and one triangular. So the plans describe 
the fundamental polyhedrons completely.  

\subsection{Describing the modular complexes for $d=1,3$} 
Let  $K$ be a Euclidean imaginary quadratic field. Then the two properties 
of the Bianchi tesselation   in the subsection 4.2 just mean that 
the groups ${\rm M}_{{\cal O}_K}^0$ and  ${\rm M}_{{\cal O}_K}^2$ are   
right $GL_2({\cal O}_K)$-modules with one generator. 
A generator can be picked as follows: 

For ${\rm M}_{{\cal O}_K}^0$: the oriented polyhedron  ${\Bbb B}_d$.

For ${\rm M}_{{\cal O}_K}^2$: the oriented geodesic ${\Bbb G} = (0, \infty)$. 

\noindent
The  $GL_2({\cal O}_K)$-module ${\rm M}_{{\cal O}_K}^1$ has one generator 
 if and only if $d=1,3$. In the latter case we have

${\rm M}_{{\cal O}_K}^1$: is generated by 
the oriented basic triangle ${\Bbb T} = (0, 1, \infty)$.

\vskip 3mm
Let  ${\Bbb D}^k_{{\cal O}_K}$ be the stabilizer of 
the corresponding generator of ${\rm M}_{{\cal O}_K}^k$ up to a sign. 
These groups are described as follows:

(i) For any imaginary quadratic number field $K$ we have 
$${\Bbb D}^2_{{\cal O}_K}
:= \mbox{the stabilizer of the geodesic 
${\Bbb G}$} \stackrel{\sim}{=} \mbox{the semidirect product of $S_2$ and ${\cal O}^*_{K}\times {\cal O}^*_{K}$}.  
$$where $S_2$ acts by permuting 
the factors. The subgroup 
${\Bbb D}^2_{{\cal O}_K} \subset GL_2({\cal O}_K)$  consists of transformations 
$$
(\sigma; \varepsilon_1, \varepsilon_2): 
 (v_1, v_2) \lms 
(\varepsilon_1 v_{\sigma(1)}, \varepsilon_2 v_{\sigma(2)}); \qquad 
\varepsilon_i \in {\cal O}^*_{K}. \quad \sigma \in S_2.
$$

(ii) For any imaginary quadratic number field $K$ there is an isomorphism
$$
{\Bbb D}^1_{{\cal O}_K} := \mbox{the stabilizer of the geodesic 
triangle ${\Bbb T}$} \stackrel{\sim}{=} S_3 \times {\cal O}^*_{K}.
$$
It can be described as follows. Let 
$L_{2}$ be a free rank two ${\cal O}_K$-module.

Choose a basis $(v_1, v_2)$ of $L_2$. Then $GL_2({\cal O}_K)$ is 
the automorphism group of $L_2$ written in this basis. We
 define a vector $v_0$ of $L_2$ 
by the formula 
$v_0+ v_1+ v_2=0$. Let $S_3$ be the group of 
 permutations of the set $\{0,1,2\}$. 
The subgroup 
${\Bbb D}^1_{{\cal O}_K}$ consists of transformations 
$$
(v_1, v_2) \lms (\varepsilon v_{\sigma(1)}, 
\varepsilon v_{\sigma(2)}), \quad \varepsilon \in {\cal O}^*_{K}, \quad \sigma \in S_3.
$$

(iii) We define the group ${\Bbb D}^0_{{\cal O}_K}$ only 
for the Euclidean fields. It is the stabilizer in 
$GL_2({\cal O}_K)$ of the basic polyhedron ${\Bbb B}_d$. 
Here is its description for the Gaussian and Eisenstein integers: 

${\Bbb D}^0_{\Z[i]} \stackrel{\sim}{=}$  
(the symmetry group of the octahedron) $\times$ ${\cal O}^*_{\Z[i]}$. 

${\Bbb D}^0_{\Z[\rho]} \stackrel{\sim}{=}$  
(the symmetry group of the tetrahedron) $\times$ ${\cal O}^*_{\Z[\rho]}$.

\vskip 3mm
The group  ${\Bbb D}^k_{{\cal O}_K}$ stabilizes 
the corresponding generator of ${\rm M}_{{\cal O}_K}^k$ up to a sign. 
The sign provides a homomorphism
$
\chi_k: {\Bbb D}^k_{d} \lra \Z
$. Summarizing, we got
\begin{proposition} \label{11.2.04.1} Let $K = \Q(\sqrt{-d})$ and $d=1,3$. 
Then we have for $k=1,2,3$:
$$
{\rm M}^k_{{\cal O}_K} = 
\Z\left[GL_2({\cal O}_K)\right]\otimes_{{\Bbb D}^k_{d}} \chi_k. 
$$
\end{proposition}
The action of the group $GL_2({\cal O}_K)$ is induced by the 
multiplication from the left.

\section{Euler complexes, CM elliptic curves, and modular hyperbolic  $3$-folds}

Let us  
restrict the Euler complex ${\Bbb E}^{\bullet}_N$ to a CM point of $X(N)$ 
corresponding to a CM curve $E_K$ with complex 
multiplication by ${\cal O}_K$, and then 
specialize further, by considering only 
${\cal N}$-torsion points.  

\subsection{Euler complexes for CM elliptic curves 
and geometry of modular orbifolds} 
{\it Elliptic units for a CM elliptic curve.} 
Let $K$ be an imaginary quadratic field, and  
 $E_{K}$ a CM elliptic curve 
with the endomorphism ring ${\cal O}_K$. Let us assume that the class number 
of $K$ is one. Then the curve $E_K$ is defined over $K$. 
The set of complex points of $E_K$ is given by $\C/\sigma({\cal O}_K)$, where 
$\sigma: K \hookrightarrow  \C$ is an embedding.  

Recall that 
${\rm Aut}(E_K) = {\cal O}_{K}^*$.  
We denote by $\varepsilon(a)$ 
the image of a point $a$ under the action of the 
automorphism of $E_K$ corresponding to $\varepsilon \in {\cal O}_{K}^*$. 
Let ${\cal N}$ be an ideal of ${\cal O}_K$. The group 
${\cal O}_{K}^*$ acts by automorphisms of the group of 
${\cal N}$-torsion points of $E_K$. 

The following result is well known, see \cite{KL}. 

\begin{lemma} \label{11.7.04.10qt} 
Let $K$ be an imaginary quadratic field, 
 $E_{K}$ the corresponding CM elliptic curve, and ${\cal N}$  an ideal in 
${\cal O}_K$. 
Then the  only relations between the elliptic units corresponding to 
 ${\cal N}$-torsion points of $E_K$ are the 
distribution relations (\ref{11.8.04.4}) and  the 
symmetry relations  provided by Corollary \ref{11.3.04.10}:
\begin{equation} \label{11.8.04.2} 
\theta_E(\varepsilon(a)) = \theta_E(a), \qquad \varepsilon \in {\cal O}_{K}^*.
\end{equation}

In particular 
if ${\cal N}$ is a prime ideal in ${\cal O}_{K}$,  
the  only relations between the corresponding elliptic units are the 
symmetry relations (\ref{11.8.04.2}). 
\end{lemma}

\vskip 3mm
Let $K_{\cal N}$ be the field obtained by adjoining to $K$ the $N$-th roots of unity and elliptic units 
$\theta_E(a)$, where $a$ runs through all nonzero ${\cal N}$-torsion points 
of the curve $E_{K}$. The field $K_{\cal N}$ 
coincides with the ray class field corresponding to the ideal ${\cal N}$, 
see \cite{KL}, Chapter 11.

\begin{definition} \label{11.6.04.1} Set 
$$
{\cal C}^{\rm un}_{1}({\cal N}) := \Bigl(\mbox{The group of elliptic units in 
${\cal O}_{K_{\cal N}}$}\Bigr)\otimes \Q \stackrel{\sim}{=}
{\cal O}_{K_{\cal N}}^*\otimes \Q, 
$$
$$
{\cal C}_{1}({\cal N}) := \Bigl(\mbox{The group of elliptic 
${\cal N}$-units in 
${\cal O}_{K_{\cal N}}$}\Bigr)\otimes \Q \stackrel{\sim}{=}
{\cal O}_{K_{\cal N}}[\frac{1}{N}]^*\otimes \Q,  
$$
$$
\widehat {\cal C}_{1}({\cal N}):= 
{\cal C}_{1}({\cal N}) \oplus \Q.
$$ 
\end{definition} 
In the last formula the factor $\Q$ formally corresponds to the non-existing element $\theta_E(0)$. We abuse notation and denote the element $1 \in \Q$ 
there by $\theta_E(0)$.

\begin{proposition} \label{erer}
Let $a_1,a_2,a_3$ be three ${\cal N}$-torsion points of $E$ with 
$a_1+ a_2+a_3=0$. Then 
$\theta_E(a_1,a_2,a_3)
 \in B_2(K_{{\cal N}})_\Q$. 
\end{proposition}

{\bf Proof}. By the very definition, 
$\theta_E(a_1,a_2,a_3) \in {\cal B}_2(\overline \Q)_\Q$. 
Recall that ${\cal B}_2(\overline \Q)_\Q = 
{B}_2(\overline \Q)_\Q$. The rest follows from Lemma 
\ref{6.11.05.123}. The proposition follows.

\begin{definition} The $\Q$-vector space 
${\cal C}_2({\cal N})$ is the subspace of  
$B_2(K_{{\cal N}})_\Q$ generated by the elements $\theta_E(a_1,a_2,a_3)$ when $a_1, a_2, a_3$ run through all 
${\cal N}$-torsion points of $E$ with $a_1+a_2+a_3=0$. 
\end{definition}

It follows from (\ref{22:10}) that there are 
the following complexes, placed in degrees $[1,2]$:
\begin{equation} \label{21:58}
{\bf E}^{\bullet}({\cal N}):= \quad 
{\cal C}_{2}({\cal N})\stackrel{\delta_2}{\longrightarrow} 
\Lambda^2{\cal C}_{1}({\cal N}); \qquad 
\widehat {\bf E}^{\bullet}({\cal N}):= \quad 
{\cal C}_{2}({\cal N})\stackrel{\delta_2}{\longrightarrow} 
\Lambda^2\widehat {\cal C}_{1}({\cal N}).
\end{equation}

\begin{lemma} \label{6.14.05.2}
Let  ${\cal N}={\cal P}$ be a prime ideal. Then 
there is one more complex: 
\begin{equation} \label{21:58}
{\bf E}_{\rm un}^{\bullet}({\cal P}):= \quad 
{\cal C}_{2}({\cal P})\stackrel{\delta_2}{\longrightarrow} 
\Lambda^2{\cal C}^{\rm un}_{1}({\cal P}); 
\end{equation}
\end{lemma}

Let $\Gamma_1(  {\cal N})$ be the subgroup of $GL_2({\cal O}_K)$ 
fixing the row vector $(0,1)$ modulo ${\cal N}$. Set 
$$
{\cal Y}_1(  {\cal N}):= \Gamma_1(  {\cal N})\backslash {\cal H}_3, \qquad 
{\cal X}_1(  {\cal N}):= \Gamma_1(  {\cal N}) \backslash {\cal H}^*_3.
$$

\vskip 3mm 
{\it The Galois group ${\rm Gal}(K_{\cal P}/K)$}. 
Let ${\Bbb F}_{\cal P} = {\cal O}_K/{\cal P}$ be the residue field of 
a prime ideal ${\cal P}$. Let $\mu_{\cal P}$ be the projection of the group of units to 
${\Bbb F}^*_{\cal P}$.  
The Galois group ${\rm Gal}(K_{\cal P}/K)$ is canonically 
isomorphic to 
${\Bbb F}_{\cal P}^*/\mu_{\cal P}$. 
 Indeed, 
$K_{\cal P}$ is 
the ray class field corresponding to the ideal ${\cal P}$, and  $h_K =1$.

\vskip 3mm
{\it The diamond operators}. Let $N(G)$ 
the normalizer of a group  $G$. Then the group 
$$
{\cal D}_{\cal P}:= N({\Gamma}_1(  {\cal P}))/{\Gamma}_1(  {\cal P})
$$ 
acts on ${\cal Y}_1(  {\cal P})$ preserving the Bianchi tesselation. Hence it acts on the modular complex. The group  ${\cal D}_{\cal P}$ is naturally isomorphic to 
${\Bbb F}_{\cal P}^*/\mu_{\cal P}$. It 
consists of the diamond operators $\langle a\rangle$, where 
$a \in {\Bbb F}_{\cal P}^*/\mu_{\cal P}$. 
Therefore there are canonical isomorphisms
$$
{\cal D}_{\cal P} = {\Bbb F}_{\cal P}^*/\mu_{\cal P} = {\rm Gal}(K_{\cal P}/K).
$$

Given a prime ideal ${\cal P}$, we define a map 
$$
\Lambda^2 \widehat {\cal C}_{1}({\cal P})
 \quad \stackrel{\delta'}{\lra } 
\quad {\cal C}_{1}({\cal P}) \oplus {\cal C}_{1}({\cal P}).
$$
by setting 
\begin{equation} \label{11.8.04.13}
\delta': \theta_E(a) \wedge \theta_E(b) \lms  
\left\{ \begin{array}{ll}
0  \oplus  \theta_E(a) - \theta_E(b)
& \quad \mbox{if } a, b \not = 0 \\ 
 - \theta_E(b)\oplus  \theta_E(b)&  \quad \mbox{if } a = 0 , b \not = 0.  \end{array} \right.
\end{equation}
Then clearly  the composition $\delta'\circ \delta_2 =0$. 
Further, there is a map  ${\cal C}_1({\cal P}) \lra \Q$, $\theta_E(a) \lms 1$. 
Taking the sum of the two copies of this map, we get a map 
$\Sigma: {\cal C}_{1}({\cal P}) \oplus {\cal C}_{1}({\cal P}) \lra \Q$. 
So we get  the following complex, placed in degrees $[1,4]$:
\begin{equation} \label{6.10.00.4}
{\cal C}_{2}({\cal P}) \quad \stackrel{\delta_2}{\lra } \quad 
\Lambda^2 \widehat {\cal C}_{1}({\cal P}) 
\quad \stackrel{\delta'}{\lra } 
\quad 
{\cal C}_{1}({\cal P}) \oplus {\cal C}_{1}({\cal P}) 
\stackrel{\Sigma}{\lra} \Q.
\end{equation}

\begin{lemma} \label{6.14.05.21}
The complex (\ref{6.10.00.4}) is quasiisomorphic to the complex 
${\bf E}_{\rm un}^{\bullet}({\cal P})$. 
\end{lemma}

{\bf Proof}. There is an obvious map of complexes
$$
\begin{array}{ccccccc}
{\cal C}_{2}({\cal P}) &\stackrel{\delta_2}{\lra }&
\Lambda^2 {\cal C}^{\rm un}_{1}({\cal P}) &&&&\\
=\downarrow && \downarrow &&&&\\
{\cal C}_{2}({\cal P}) &\stackrel{\delta_2}{\lra } &
\Lambda^2 \widehat {\cal C}_{1}({\cal P}) 
& \stackrel{\delta'}{\lra } 
&
{\cal C}_{1}({\cal P}) \oplus {\cal C}_{1}({\cal P})&
\stackrel{\Sigma}{\lra } & \Q.
\end{array}
$$
It provides the desired quasiisomorphism. Indeed, 
it follows from the very definition that 
 $\widehat {\cal C}_{1}({\cal P}) = {\cal C}^{\rm un}_{1}({\cal P}) 
\oplus \Q \oplus \Q$. It is easy to see from this and the definition 
of the map $\delta'$ that  bottom complex is exact in degree $3$, 
and the kernel of the map $\delta'$ is 
$\Lambda^2{\cal C}^{\rm un}_{1}({\cal P})$. 
The lemma is proved. 

\begin{theorem} \label{mtheor2da} Let us assume that $K$ is either 
Gaussian or Eisenstein field, i.e. $d=1,3$. Let ${\cal P}$ be 
a prime ideal in ${\cal O}_K$. Choose a ${\cal P}$-torsion 
point $z$ of the curve $E_K$ generating $E_K[{\cal P}]$. 

Then there exists a surjective homomorphism  of complexes
\begin{equation} \label{6.16.05.10}
\begin{array}{ccccc}
\Q\otimes _{\Gamma_1(  
{\cal P})}\Bigl({\rm M}^1_{{\cal O}_K}&\stackrel{\partial}{\longrightarrow}& 
{\rm M}^2_{{\cal O}_K}&\stackrel{\partial}{\longrightarrow}&{\rm M}^3_{{\cal O}_K}\Bigr)\\
&&&&\\
\downarrow\theta^{(1)}  &&\downarrow\theta^{(2)} & &\downarrow\theta^{(3)} \\
&&&&\\
{\cal C}_{2}({\cal P})&\stackrel{\delta_2}{\longrightarrow}& 
\Lambda^2\widehat {\cal C}_{1}({\cal P})
&\stackrel{\delta'}{\lra} &{\cal C}_{1}({\cal P}) \oplus {\cal C}_{1}({\cal P})
\end{array}
\end{equation}
which intertwines the action of the group  ${\cal D}_{\cal P}$
with the action of ${\rm Gal}(K_{\cal P}/K)$. 

Moreover the maps $\theta^{(2)}$ and $\theta^{(3)}$ are isomorphisms. 
\end{theorem}

{\bf Proof}. Since ${\cal O}_K$ has the class number one,  one has 
$$
\Gamma_1(  {\cal P})\backslash GL_2({\cal O}_K) = {\Bbb F}_{\cal P}^2 - \{0,0\} = 
\{(\alpha, \beta) \in {\Bbb F}_{\cal P}^2- 0\}.
$$ 
We identify it with the set 
of nonzero rows $\{(\alpha, \beta, \gamma) \}$ 
with $\alpha + \beta + \gamma = 0$. 
Then 
\begin{equation} \label{11.48.04.7}
\Q\otimes_{\Gamma_1(  
{\cal P})}{\rm M}_{{\cal O}_K}^1=  \Z[\Gamma_1(  {\cal P})\backslash GL_2({\cal O}_K)]
\otimes_{{\Bbb D}_{{\cal O}_K}^1} \chi_1  =
\end{equation}
$$
\frac{\Z[(\alpha, \beta, \gamma)\in {\Bbb F}_{\cal P}^3- 0 \quad | \quad \alpha + 
\beta + \gamma = 0]}{\mbox{the dihedral symmetry relations } }.
$$
where the dihedral symmetry relations are the following:
$$
 (\alpha_1, \alpha_2, \alpha_3) = (-1)^
{{\rm sgn}(\sigma)}(\alpha_{\sigma(1)}, 
\alpha_{\sigma(2)}, \alpha_{\sigma(3)}), \qquad \sigma \in S_3, 
$$
$$
 (\alpha_1, \alpha_2, \alpha_3) = (\varepsilon \alpha_{1}, \varepsilon \alpha_{2}, \varepsilon \alpha_{3}), \qquad \varepsilon \in {\cal O}_{K}^*.
$$
\vskip 3mm 
\begin{equation} \label{11.48.04.8}
\Q\otimes_{\Gamma_1(
{\cal P})}{\rm M}_{{\cal O}_K}^2=  \Z[\Gamma_1(  {\cal P})\backslash GL_2({\cal O}_K))]
\otimes_{{\Bbb D}_{{{\cal O}_K}}^2} \chi_2 =
\end{equation}
$$ \frac{\Z[ (\alpha, \beta)\in {\Bbb F}_{\cal P}^2-0]}{ (\alpha, \beta) =  
-(\beta, \alpha) = (\varepsilon \alpha, \varepsilon \beta) }, 
\qquad \varepsilon \in {\cal O}_K^*.
$$

\vskip 3mm 

\begin{equation} \label{11.48.04.9}
\Q\otimes_{\Gamma_1(2, 
{\cal P})}{\rm M}_{{\cal O}_K}^3 = 
\Z[\Gamma_1( {\cal P})\backslash GL_2({\cal O}_K)/ 
B({\cal O}_K)]  = 
\end{equation}
$$
\mbox{the free abelian group with generators $[\beta, 0]$ and $[0, \beta]$,  where 
$\beta \in {\Bbb F}_{\cal P}^*/\mu_{\cal P}$}.
$$
Here $B $ is group of upper triangular matrices. 
\vskip 3mm

Then the  differential in the top complex in (\ref{11.48.04.6}) 
is described  as follows:
\begin{equation} \label{11.48.04.10}
\partial: (\alpha, \beta, \gamma) \lms (\alpha, \beta) + 
(\beta, \gamma) + (\gamma, \alpha); 
\end{equation}
\begin{equation} \label{11.48.04.11}
\partial: (\alpha, \beta)\lms  
\left\{ \begin{array}{ll}
0 \quad \oplus  [0, \alpha] - [0, \beta]
& \quad \mbox{if } \alpha, \beta \not = 0 \\ 
 - [\beta, 0]\quad \oplus  \quad  [0, \beta]&  \quad \mbox{if } \alpha = 0 , \beta \not = 0.  \end{array} \right.
\end{equation}
\vskip 3mm

To define the 
map  of complexes $\theta^{(*)}$ we pick a nonzero ${\cal P}$-torsion 
point $z$ of the curve $E_K$. Then there is an isomorphism of abelian groups
\begin{equation} \label{11.7.04.3}
{\Bbb F}_{\cal P} \lra E_K[{\cal P}], \qquad \alpha \lms \alpha z.
\end{equation}
Now we proceed as follows.

(i) Let ${\cal F}_{\cal P}:= \Q\left[{\Bbb F}^*_{\cal P}/\mu_{\cal P}\right]$ be the $\Q$-vector space generated by the set 
${\Bbb F}^*_{\cal P}/\mu_{\cal P}$. 
We denote by $[\alpha]$ the generator of ${\cal F}_{\cal P}$ 
corresponding to an 
element $\alpha \in {\Bbb F}^*_{\cal P}$, so $[\varepsilon\alpha] = [\alpha]$. 
According to Proposition \ref{11.6.04.1} there is an isomorphism
\begin{equation} \label{11.7.04.2}
{\cal F}_{\cal P}\stackrel{\sim}{\lra} \Bigl(\mbox{The group of 
${\cal P}$-elliptic units in 
${\cal O}^*_{K_{\cal P}}$}\Bigr)\otimes \Q,  \qquad 
  [\alpha]\lms\theta_E(\alpha z).
\end{equation}

(ii) Using the identifications (\ref{11.48.04.7}) - (\ref{11.48.04.9}) 
we define the desired 
map  of complexes as follows: 
\begin{equation} \label{6.17.05.11}
\theta^{(1)}: (\alpha, \beta, \gamma) \lms  
\theta_E(\alpha z, \beta z, \gamma z),
\end{equation}
\begin{equation} \label{6.17.05.12}
\theta^{(2)}: (\alpha, \beta) \lms  
\theta_E(\alpha z)\wedge \theta_E(\beta z),
\end{equation}
$$
\theta^{(3)}: [\beta_1, 0] + [0, \beta_2] \lms \theta_E(\beta_1 z) \oplus 
\theta_E(\beta_2 z).
$$
The map $\theta^{(1)}$ is well defined thanks to the dihedral symmetry relations 
for $\theta_E(a, b, c)$; the maps $\theta^{(2)}$ and $\theta^{(3)}$ are
 well defined thanks to the 
 symmetry relations 
for $\theta_E(a)$. 
It follows easily from (\ref{11.48.04.10})-(\ref{11.48.04.11}), the formula 
(\ref{22:10}),  
and (\ref{11.8.04.13}) 
that the map $\theta^{(*)}$ is 
a homomorphism of 
complexes. It is surjective by Theorem \ref{21:56} b) and Proposition 
\ref{11.6.04.1}. 
Since ${\cal P}$ is a 
prime ideal, thanks to Lemma
\ref{11.7.04.10qt} 
there are no distribution relations between the elliptic units, 
and hence the maps 
$\theta^{(2)}$ and $\theta^{(3)}$ are isomorphisms.

It follows from the very definition that the map $\theta^*$ intertwines the actions of the groups ${\cal D}_{\cal P}$ and ${\rm Gal}(
K_{\cal P}/K)$. 
The theorem is proved. 
\vskip 3mm 

{\bf Remark}. Let us 
remove dependence on the choice of a ${\cal P}$-torsion point $z$ in the 
theorem. Let $\mu(1)$ be   the ${\rm Gal}(K_{\cal P}/K)$-module 
of non-trivial ${\cal P}$-torsion points on the curve $E_K$. 
Then there is a {\it canonical} morphism 
$$
\mbox{the top complex in (\ref{6.16.05.10})} \to 
\mbox{the bottom complex in (\ref{6.16.05.10})}  \otimes_{{\rm Gal}(K_{\cal P}/K)}\mu(1).
$$
 It is given by 
$\theta^*_z \otimes z^{-1}$, where $\theta^*_z$ 
is the above map of complexes defined using a torsion point $z$. 

\subsection{The rational cohomology of $\Gamma_1(  {\cal P})$} 
If $\Gamma$ is a torsion free finite index subgroup of $GL_2({\cal O}_K)$, 
we have an isomorphism
$$
H^*(\Gamma, \Q) = H^*(\Gamma \backslash {\cal H}_3, \Q). 
$$
If we understood the right hand side as the 
cohomology group of the orbifold 
${\cal Y}_{\Gamma}:= \Gamma \backslash {\cal H}_3$, 
this  formula remains valid for any subgroup of $GL_2({\cal O}_K)$. 
Let us recall few basic facts about the cohomology groups 
$H^*(\Gamma, \Q)$. 
See \cite{Ha1} and references there for further details. 
 
The orbifold ${\cal Y}_{\Gamma}$ is compactified  by 
${\cal X}_{\Gamma}:= \Gamma \backslash {\cal H}^*_3$. 
The complement  ${\cal X}_{\Gamma} - {\cal Y}_{\Gamma}$ 
consists of a finite number of cusps. The orbifold 
${\cal Y}_{\Gamma}$ has a boundary 
$\partial {\cal Y}_{\Gamma}$  given by disjoint union 
of two dimensional orbifold tori. These tori are parametrized by the cusps. 
Restriction to the boundary provides a map
\begin{equation} \label{mm}
{\rm Res}: H^*({\cal Y}_{\Gamma}, \Q) \lra H^*(\partial {\cal Y}_{\Gamma}, \Q).
\end{equation}
The kernel of the restriction map is the cuspidal part 
$H_{\rm cusp}^*(\Gamma, \Q)$ of the cohomology. The image is called 
the Eisenstein part of the cohomology, and denoted 
$H_{\rm Eis}^*(\Gamma, \Q)$. So 
there is  an exact sequence
$$
0 \lra H_{\rm cusp}^*(\Gamma, \Q) \lra 
H^*(\Gamma, \Q) \lra 
H_{\rm Eis}^*(\Gamma, \Q) \lra 0.
$$
The cuspidal cohomology satisfies the Poincare duality, so we have 
\begin{equation} \label{11.11.04.1} 
{\rm dim} H_{\rm cusp}^1(\Gamma, \Q) = 
{\rm dim} H_{\rm cusp}^2(\Gamma, \Q).
\end{equation}

\vskip 3mm
Let $M$ be a threefold with boundary $\partial M$. Then 
it follows from the Poincare duality that 
the image of the restriction map ${\rm Res}: H^*(M) \to H^*(\partial M)$ 
is a Lagrangian subspace 
in $H^*(\partial M)$. 

Therefore if $\Gamma$ is torsion free then 
the image of the restriction map is a Lagrangian subspace in 
the boundary cohomology
$H^*(\partial {\cal Y}_{\Gamma}, \Q)$. 

It follows from this that if $\Gamma$ is torsion free then 
the homological Euler characteristic is zero:
$$
\chi_h(\Gamma) = \sum(-1)^i{\rm dim}H^i(\Gamma, \Q) = 0.
$$

Indeed, clearly $H^0_{\rm Eis}(\Gamma, \Q) =\Q$. Since 
the image of the restriction map (\ref{mm}) is Lagrangian, this implies that 
there is a natural isomorphism
$$
H_{\rm Eis}^2(\Gamma, \Q) = {\rm Ker}
\Bigl(\Q[\mbox{cusps of $\Gamma$}] \stackrel{\Sigma}{\lra} \Q\Bigr).
$$
Observe that identifying $H_{\rm Eis}^2(\Gamma, \Q)$ with the Borel-Moore 
homology $H_1^{BM}({\cal Y}_{\Gamma}, \Q)$, the boundary map is 
given by the ``boundary at infinity''
 of  the $1$-cycle representing the class, which is a 
linear combination of cusps of total degree zero. 

Finally, ${\rm dim}H^1(\partial {\cal Y}_{\Gamma}) = 
2 \times \mbox{(the number of cusps of $\Gamma$)}$, and therefore 
$
{\rm dim} H_{\rm Eis}^1(\Gamma, \Q)$ is equal to the number of cusps 
of $\Gamma$. 
However if $\Gamma$ has torsion the situation is more complicated.

\vskip 3mm
{\bf Example}. The modular $3$-fold ${\cal Y}_1({\cal P})$ is a three 
dimensional orbifold. Indeed, 
the subgroup of diagonal matrices ${\rm diag}(\varepsilon, 1)$, where 
$\varepsilon \in {\cal O}_K^*$,  
is a torsion subgroup in $\Gamma_1(  {\cal P})$. 
The complement ${\cal X}_1({\cal P}) - {\cal Y}_1({\cal P})$ 
consists of $2|{\Bbb F}_{\cal P}^*/\mu_{\cal P}|$ cusps. 
The natural covering ${\cal X}_1({\cal P})\lra {\cal X}_0({\cal P})$ is 
unramified of degree $|{\Bbb F}_{\cal P}^*/\mu_{\cal P}|$. 
There are two cusps on ${\cal X}_0({\cal P})$: 
the $0$ and $\infty$ cusps. So there are $|{\Bbb F}_{\cal P}^*/\mu_{\cal P}|$ 
cusps over $0$ as well as   over $\infty$. 

\begin{lemma} \label{11.18.04.1}
Let ${\cal P}$ be a prime ideal in ${\cal O}_K$. Then 
$$
H_{\rm Eis}^1(\Gamma_1(  {\cal P}), \Q) =0,
$$
$$
H_{\rm Eis}^2(\Gamma_1(  {\cal P}), \Q) = {\rm Ker}
\Bigl(\Q[\mbox{{\rm cusps of} $\Gamma_1(  
{\cal P})$}] \stackrel{\Sigma}{\lra} \Q\Bigr).
$$
\end{lemma}

{\bf Proof}. In this case the boundary cohomology group 
$H^1$ is zero. Indeed, the cusps, and hence the boundary components 
 are parametrized by the 
$\Gamma_1(  {\cal P})$-orbits on the set of Borel subgroups in 
$GL_2({\cal O}_K)$. Consider the boundary component 
corresponding to the standard (upper triangular) 
Borel subgroup of $GL_2({\cal O}_K)$. We identify $H_1$ of the 
corresponding boundary component with the unipotent radical 
$U({\cal O}_K) = {\cal O}_K$. 
The torsion element 
${\rm diag}(\varepsilon, 1)$ acts on it by multiplication on
 $\varepsilon$. Taking  $\varepsilon =-1$, we see that 
it kills the boundary $H^1$. 
However on $H^2$ it acts as the identity. The lemma follows. 

\vskip 3mm
{\bf Remark}. This lemma is consistent with Horozov's computation of 
the homological Euler characteristic of $\Gamma_1(  {\cal P})$, see 
Theorem 0.5 in \cite{H}. Namely, 
\begin{equation} \label{345}
{\rm dim}H_{\rm Eis}^2(\Gamma_1(  {\cal P}), \Q) = 
2|{\Bbb F}_{\cal P}^*/\mu_{\cal P}|-1, 
\end{equation}
and thus we get 
$$
\chi_h(\Gamma_1(  {\cal P})) = 
1-
{\rm dim}H_{\rm Eis}^1(\Gamma_1(  {\cal P}), \Q) + 
{\rm dim}H_{\rm Eis}^2(\Gamma_1(  {\cal P}), \Q) = 
2|{\Bbb F}_{\cal P}^*/\mu_{\cal P}|
$$
which coincides with the result obtained in loc. cit. 

\subsection{Applications} 
\begin{corollary} \label{11.7.04.20} Let $K = \Q(\sqrt{-d})$ where $d=1,3$, and 
let ${\cal P}$ be a prime ideal in ${\cal O}_K$. Then there are canonical isomorphisms 
\begin{equation} \label{11.7.04.21} 
\widehat \theta_h^{(2)}: H^2(\Gamma_1(  {\cal P}), \Q) 
\stackrel{\sim}{\lra} H^2(\widehat {\bf E}^{\bullet}({\cal P})), 
\end{equation}
\begin{equation} \label{11.7.04.21a} 
\theta_h^{(2)}: H_{\rm cusp}^2(\Gamma_1(  {\cal P}), \Q) 
\stackrel{\sim}{\lra} H^2({\bf E}_{\rm un}^{\bullet}({\cal P})). 
\end{equation}
\end{corollary}

{\bf Proof}. The isomorphism (\ref{11.7.04.21}) 
follows immediately from Theorem \ref{mtheor2da}. 
Indeed, take the stupid truncation in degrees $\leq 2$ of the two complexes in 
(\ref{6.16.05.10}).
Then $H^2$ of the truncated top complex is 
$H^2(\Gamma_1(  {\cal P}), \Q)$, while  $H^2$ of the truncated bottom complex 
is $H^2(\widehat {\bf E}^{\bullet}({\cal P}))$. 

To  prove (\ref{11.7.04.21a}) observe that $H^2$ of the top complex 
delivers the left hand side in (\ref{11.7.04.21a}). 
By Lemma
\ref{6.14.05.21},  the 
$H^2$ of the bottom complex 
gives the right hand side. 

Here is another argument for this. 
There is a natural map
$$
H^2(\widehat {\bf E}^{\bullet}({\cal P})) \lra {\rm Ker}(
{\cal C}_1{\cal P} \oplus {\cal C}_1{\cal P} \stackrel{\Sigma}{\lra}\Q).
$$
The right hand side can be identified with 
$H^2_{\rm Eis}(\Gamma_1({\cal P},\Q)$, see Section 5.2, e.g. 
(\ref{345}), 
and thus the kernel is identified with $H^2_{\rm cusp}(\Gamma_1({\cal P}),\Q)$. 
 The corollary is proved.

\vskip 3mm

The commutativity of the diagram (\ref{6.16.05.10}) implies 
that the composition 
$$
 {\rm M}^0_{{\cal O}_K} \stackrel{\partial}{\lra} {\rm M}^1_{{\cal O}_K} 
\stackrel{\theta^{(1)}}{\lra} {\cal C}_{2}({\cal P})
$$
 provides 
a natural map 
\begin{equation} \label{21:52}
{\rm M}^0_{{\cal O}_K} \lra {\rm Ker}\Bigl(
{\cal C}_{2}({\cal P}) \stackrel{\delta_2}{\lra}\Lambda^2
 {\cal C}_{1}({\cal P})\Bigr) = 
H^1({\bf E}^{\bullet}({\cal P})). 
\end{equation}

\begin{definition} \label{21:51}
$\alpha(K_{\cal P})$ is the image of the map (\ref{21:52}). 
\end{definition}

Since $K_{\cal P}$ is a number field, 
${\cal B}_2(K_{\cal P})_\Q = {B}_2(K_{\cal P})_\Q$. 
So ${\rm Ker}\delta_2 =  K_3^{ }(K_{\cal P})_{\Q}$ 
by Suslin's theorem. Thus 
there is an embedding
$$
H^1({\bf E}^{\bullet}({\cal P})) \subset
K_3^{ }(K_{\cal P})_{\Q}.
$$
So we immediately get 

\begin{corollary} \label{21:59}
Under the same assumptions as in Corollary \ref{11.7.04.20}, 
there is a homomorphism
\begin{equation} \label{11.7.04.22}
\theta_h^{(1)}: H_{\rm cusp}^1(\Gamma_1(  {\cal P}), \Q) 
\stackrel{}{\lra} \frac{H^1({\bf E}^{\bullet}({\cal P})) }{
\alpha(K_{\cal P})}\subset \frac{K_3^{ }(K_{\cal P})_\Q}
{\alpha(K_{\cal P})}.
\end{equation}
\end{corollary}

{\bf Proof}. Indeed, there is a map 
\begin{equation} \label{11.7.04.22a}
\theta_h^{(1)}: H^1(\Gamma_1(  {\cal P}), \Q) 
\stackrel{}{\lra} \frac{H^1({\bf E}^{\bullet}({\cal P})) 
}{\alpha(K_{\cal P})}.
\end{equation}
Since $H_{\rm Eis}^1(\Gamma_1(  {\cal P}), \Q) = 0$, one has 
$H_{\rm cusp}^1 = H^1$, so the map
 (\ref{11.7.04.22a}) coincides with the needed map (\ref{11.7.04.22}). 
\vskip 3mm 

\subsection{Relations between the elements $\theta_E(a_1,a_2,a_3)$} 
To describe the vector space ${\cal C}_{2}({\cal N})$ 
one needs to know all relations between the elements 
$\theta_E(a_1, a_2, a_3)$. 

There is a decomposition 
\begin{equation} \label{21:53}
{\cal Y}_1({\cal P}) = \cup \overline {\gamma {\Bbb B}_d}, \qquad 
\gamma \in \Gamma_1(  {\cal P})\backslash 
PGL_2({\cal O}_K)
\end{equation}
of the  orbifold ${\cal Y}_1({\cal P})$ into a union of finite number of 
geodesic polyhedrons, obtained by projecting the Bianchi 
tessellation of ${\cal H}_3$ corresponding to  ${\cal O}_K$ 
onto  ${\cal Y}_1({\cal P})$.  
We call it the Bianchi decomposition. 
Each polyhedron of the Bianchi decomposition of ${\cal Y}_1({\cal P})$ 
gives rise to a relation 
between the elements $\theta_E(a_1, a_2,  a_3)$, where $a_i$ are 
${\cal P}$-torsion points of $E_K$. Namely, the sum of the elements 
corresponding via the map $\theta^{(1)}$  
to the faces of the polyhedron $\gamma {\Bbb B}_d$ is an element 
$$
\theta_E^{(1)}(\partial (\gamma {\Bbb B}_d)) \in 
K_3^{ }(K_{\cal P})_\Q.
$$

 We call the quotient of the 
space of relations between the elements $\theta_E(a_1,a_2,a_3)$, where $a_1,a_2,a_3$ run through the nonzero 
${\cal P}$-torsion points in $E_K$, modulo  the subspace generated by 
the dihedral symmetry relations and the described above 
relations corresponding to the polyhedrons of the decomposition 
(\ref{21:53}), the space of 
{\it sporadic level ${\cal P}$ relations} between the elements $\theta_E(a_1,a_2,a_3)$.  

\begin{corollary} \label{11.8.04.1} Let ${\cal P}$ be a prime ideal 
in ${\cal O}_K$. Then the space of level ${\cal P}$ 
sporadic relations between the elements $\theta_E(a_1,a_2,a_3)$ 
is identified with the kernel of the 
map $\theta_h^{(1)}$. 
\end{corollary}

To describe the relation corresponding to a polyhedron of the 
decomposition (\ref{21:53}) precisely 
one needs to know explicitly the element 
of $K_3^{ }(K_{\cal P})$ corresponding to this polyhedron.

\vskip 3mm 
 
\begin{conjecture} \label{11.7.04.25}
The element $\theta_E^{(1)}(\partial (\gamma {\Bbb B}_d))$ is zero.
\end{conjecture}


Summarizing, there are three types of the relations between the elements 
$\theta_E(a_1,a_2,a_3)$:
\vskip 3mm
(i) The dihedral symmetry relations (\ref{22:101}) - (\ref{22:102}) 
corresponding to 
symmetries of the geodesic triangles in the Bianchi 
decomposition of the modular orbifold ${\cal Y}_1({\cal P})$.  

(ii) The relations provided by the polyhedrons 
of the Bianchi decomposition (\ref{21:53}).

(iii) Sporadic relations are described by the map 
(\ref{11.7.04.22}), that is,    
\begin{equation} \label{11.7.04.22d}
\theta_h^{(1)}: H_{\rm cusp}^1(\Gamma_1(  {\cal P}), \Q) 
\stackrel{}{\lra} K_3^{ }(K_{\cal P})_\Q/\alpha(K_{\cal P}).
\end{equation}
\vskip 3mm 
\noindent 
 The Hecke algebra acts on the left. 
It is hard to imaging an action of the Hecke algebra on the right. 

\begin{problem} \label{11.20.04.1} Determine the map 
(\ref{11.7.04.22d}). 
\end{problem} 

\vskip 3mm
{\it Appendix}. Here is another way to assign to translations of the 
basic polyhedron elements of $K_3\otimes \Q$; this time 
they lie in the one dimensional space $K_3(K)_\Q$. 

\begin{lemma} \label{5.15} For any $\gamma \in PGL_2({\cal O}_K)$, 
the geodesic polyhedron $\gamma {\Bbb B}_d$ 
provides an element 
\begin{equation} \label{11.7.04.26}
\left[\gamma{\Bbb B}_d\right] \in K_3(K)_\Q.
\end{equation}
The value of the regulator map 
$K_3(K) \to \R$ 
on this element is the volume of the polyhedron ${\Bbb B}_d$.  
Thus the element (\ref{11.7.04.26}) does not depend on $\gamma$. 
\end{lemma}

 {\bf Proof}. Cut the polyhedron 
$\gamma {\Bbb B}_d$ into ideal geodesic simplices. 
We assine to an ideal geodesic simplex $I(\infty, 0, 1, z)$ 
with vertices at 
$\infty, 0, 1, z$, where $z \in \overline \Q$, 
the element $\{z\}_2$ of the Bloch group $B_2(\overline \Q)$. Let 
$\left[\gamma{\Bbb B}_d\right]$  be the sum 
of the elements corresponding to the simplices of the decomposition. 
It lies in $B_2(K)$ since the vertices of the simplices can be taken 
at the cusps, which are identified with  $K\cup \{\infty\}$. 
Let us show that the differential $\delta_2$ in the Bloch complex 
kills this element.

The Dehn invariant of the ideal geodesic simplex $I(\infty, 0, 1, z)$ 
equals  
$$
\log |1-z| \otimes {\rm arg} (z) - \log |z| \otimes {\rm arg} (1-z) 
\in \R \otimes S^1 = (\Lambda^2\C^*)^-
$$ 
where $-$ stands for the antiinvariants of the complex conjugation acting on 
$\Lambda^2\C^*$. A union of finite number of translations of the polyhedron 
${\Bbb B}_d$ gives (a multiple of) the fundamental cycle of the 
threefold ${\cal Y}_1({\cal P})$.  The Dehn invariant of the latter is zero. 
Thus the 
Dehn invariant of ${\Bbb B}_d$ 
 is torsion. 

According to the Lobachevsky formula one has 
$$
{\rm vol}I(\infty, 0, 1, z) = {\cal L}_2(z).
$$
Since the regulator map on $K_3^{ }(K)_\Q$ is given by the 
dilogarithm function ${\cal L}_2(z)$, 
we get the second claim of the Lemma. 
Now the injectivity of the regulator map on $K_3^{ }(K)_\Q$ implies that 
the element (\ref{11.7.04.26}) does not
depend on  $\gamma$. The Lemma is proved.

\vskip 3mm
{\bf Remark}. The element $[\gamma {\Bbb B}_d]$ can not 
be a non-zero multiple of $\theta_E^{(1)}(\partial ({\Bbb B}_d))$. 
Indeed, consider the sum of (translations of)  
basic polyhedrons providing the fundamental cycle of 
the threefold ${\cal Y}_1({\cal P})$. Then the sum of the 
corresponding elements $\theta_E^{(1)}(\partial (\gamma {\Bbb B}_d))$
 is zero. Indeed, the contribution of each 
polyhedron is provided by the triangles forming its boundary, 
and the boundary of the fundamental cycle is zero. On the other hand 
the sum of the elements $[\gamma {\Bbb B}_d]$ corresponding to 
the polyhedrons entering to the fundamental cycle is, by Lemma \ref{5.15}, a non-zero multiple of $[{\Bbb B}_d]$.

\subsection{Generalizations}  
\begin{theorem} \label{mtheor2d} Let ${\cal O}_K$ be the ring of 
Gaussian or Eisenstein integers. Let ${\cal N}$ be an ideal in ${\cal O}_K$. Then there exists a canonical surjective homomorphism  of complexes
\begin{equation} \label{11.48.04.6}
\begin{array}{ccc}
\Q\otimes _{\Gamma_1(  
{\cal N})}\Bigl({\rm M}^1_{{\cal O}_K}&\stackrel{\partial}{\longrightarrow}& 
{\rm M}^2_{{\cal O}_K}\Bigr)\\
\downarrow\theta^{(1)}  &&\downarrow\theta^{(2)}  \\
{\cal C}_{2}({\cal N})&\stackrel{\delta_2}{\longrightarrow}& 
\Lambda^2\widehat {\cal C}_{1}({\cal N})
\end{array}
\end{equation}
The maps $\theta^{(*)}$ intertwine the action of the group  ${\cal D}_{\cal N}$
with the action of ${\rm Gal}(K_{\cal N}/K)$. 
\end{theorem}

{\bf Proof}. We follow the proof of Theorem \ref{mtheor2da}
 with necessary modifications. One has 
$$
\Gamma_1(  {\cal N})\backslash GL_2({\cal O}_K) = 
\{(\alpha, \beta) \in ({\cal O}_K/{\cal N})^2 | (\alpha, \beta, {\cal N})=1\}.
$$ 
We identify it with the set 
of rows $\{(\alpha, \beta, \gamma) \}$ 
with $\alpha + \beta + \gamma = 0$, and $\alpha, \beta, \gamma$ are modulo 
${\cal N}$ and have no common divisors with ${\cal N}$.  
Then 
\begin{equation} \label{11.48.04.7s}
\Q\otimes_{\Gamma_1(  
{\cal N})}{\rm M}_{{\cal O}_K}^1 \quad =\quad   
\frac{\Z[(\alpha, \beta, \gamma)\in ({\cal O}_K/{\cal N})^3 
\quad | \quad \alpha + 
\beta + \gamma = 0, \quad (\alpha, \beta, {\cal N})=1]}{\mbox{the dihedral symmetry relations } },
\end{equation}
\begin{equation} \label{11.48.04.8}
\Q\otimes_{\Gamma_1( {\cal N})}{\rm M}_{{\cal O}_K}^2 \quad =\quad 
\frac{\Z[ (\alpha, \beta)\in ({\cal O}_K/{\cal N})^2] \quad | \quad 
(\alpha, \beta, {\cal N}) = 1)}{ (\alpha, \beta) =  
-(\beta, \alpha) = (\varepsilon \alpha, \varepsilon \beta) }, 
\qquad \varepsilon \in {\cal O}_K^*.
\end{equation}

To define the maps $\theta^{(*)}$, we pick a primitive ${\cal N}$-torsion point $z$ on $E_K$, and use formulas 
(\ref{6.17.05.11}) - (\ref{6.17.05.12}). 
Then we get a homomorphism of complexes. To show that it is surjective, we use the distribution relations. The theorem is proved.

\vskip 3mm

It would be very interesting to find the ``right'' 
generalisation of these results for other imaginary quadratic fields. 
The above construction has an obvious generalization: the subcomplex 
of the modular complex, given by the 
$GL_2({\cal O}_K)$-submodules generated by ${\Bbb T}$ and ${\Bbb G}$, 
is mapped to the complex ${\cal C}_{2}({\cal N})\stackrel{\delta_2}
{\longrightarrow}
\Lambda^2\widehat {\cal C}_{1}({\cal N})$. However 
the $GL_2({\cal O}_K)$-submodules generated by ${\Bbb T}$ and 
${\Bbb G}$ do not give, in general, the whole  modular complex. 
For example, let 
$K = \Q(\sqrt{-2})$. Then  
the $GL_2({\cal O}_K)$-submodule generated by 
${\Bbb G}$ is ${\cal M}^2_{{\cal O}_K}$. However 
there are $2$-cells of the Bianchi tessellation 
which are not in the $GL_2({\cal O}_K)$-orbit of ${\Bbb T}$, e.g. the geodesic $4$-gon with vertices at 
$(1, \frac{1+\theta}{2}, 1+\theta, \infty)$. 
This suggsets that there might exist 
 a natural construction of 
an element of $B_2({\cal S}_{\cal N})$, which goes under the differential $\delta_2$ to 
  the sum of four wedge products of 
the elliptic units corresponding to the sides of this $4$-gon, and similalrly for other $2$-cells 
in the $GL_2({\cal O}_K)$-orbit of this $4$-gon. 

Another question is whether/when the $1$-cycles provided by the $GL_2({\cal O}_K)$-orbits of ${\Bbb G}$ 
generated the Borel-Moore $H_1$ of the corresponding modular $3$-fold for general $K$. 
If so, this would imply that for a prime ideal ${\cal P}$ of ${\cal O}_K$ there is a 
natural surjective map from 
$\Lambda^2{\cal O}^*_{K_{\cal P}}$ to $H^2_{\rm cusp}$ of the corresponding modular $3$-fold, 
just like in Conclusion 2.

Dept of Mathematics, Brown University, Providence RI 02912, USA. e-mail: sasha@math.brown.edu

\end{document}